\documentclass{article}

\usepackage[english]{babel}

\usepackage[letterpaper,top=2cm,bottom=2cm,left=3cm,right=3cm,marginparwidth=1.75cm]{geometry}
\usepackage{enumitem}
\usepackage{amsthm}
\usepackage{subfigure}
\usepackage{caption}

\usepackage{amsmath}
\usepackage{amssymb}
\usepackage{graphicx}
\usepackage{bm}
\usepackage[colorlinks=true, allcolors=blue]{hyperref}
\usepackage[linesnumbered,ruled,vlined]{algorithm2e}
\usepackage{comment}
\usepackage{float}
\renewcommand{\div}{\text{div}}
\newcommand{\divt}{\widetilde{\text{div}}}
\newcommand{\bv}{\bm{v}}
\newcommand{\V}{\bm{V}}

\newcommand{\x}{\bm{x}}

\renewcommand{\i}{\mathrm i}
\newcommand{\Sspace}{{\mathbb S}}

\usepackage{xcolor}

 \newcommand{\degp}{d}

\newcommand{\Rr}{\mathbb{R}}
\newcommand{\Nn}{\mathbb{N}}
\newcommand{\Ran}{\mathcal{R}}
\newcommand{\Ker}{\text{Ker}}
\newcommand{\Span}{\text{Span}}

\newcommand{\al}[1]{\textcolor{blue}{#1}}

\theoremstyle{definition}
\newtheorem{definition}{Definition}[subsection]

\newtheorem{proposition}[definition]{Proposition}

\theoremstyle{remark}
\newtheorem{remark}[definition]{Remark}

\title{Quasi-Trefftz spaces for a first-order formulation of the Helmholtz equation }
\author{Lise-Marie Imbert-G\'erard, Andr\'ea Lagard\`ere, Guillaume Sylvand, S\'ebastien Tordeux}
\date{}

\begin{document}
\maketitle

\tableofcontents
\newpage

{\bf Abstract}\\
This work is concerned with the development of quasi-Trefftz methods for first-order differential systems. It focuses on discrete quasi-Trefftz spaces, starting from their definition and including the construction of corresponding bases together with their computational aspect.

This is the first attempt at constructing quasi-Trefftz bases for a problem governed by a first-order system without relying on an auxiliary scalar equation. 
A decoupling approach, with a second order scalar equation for the one unknown, is proposed here simply as a point of comparison to this new approach.

\section*{Acknowledgments}
L.-M. Imbert-G\'erard acknowledges support from the US National Science Foundation: this material is based
upon work supported by the United States National Science Foundation under Grant No.  DMS-2110407.

\section*{Statements and declarations}
L.-M. Imbert-G\'erard has disclosed an outside interest in Airbus Central R\&T to the University of Arizona.
Conflicts of interest resulting from this interest are being managed by The University of Arizona in accordance
with its policies.

\section{Introduction}

Trefftz Discontinuous Galerkin (DG) methods are numerical methods for the approximation of solutions to Partial Differential Equation (PDE) problems, see for instance \cite{Hiptmair2016}. As Galerkin methods, their specificity is to leverage spaces of local exact solutions to the governing PDE, so-called Trefftz functions, for trial and test function spaces alike, both at the continuous and discrete levels. 
There are two crucial comparison points between discrete spaces of Trefftz functions and standard polynomial spaces in terms of best-approximation properties:
\begin{itemize}
    \item while polynomial spaces can approximate smooth enough functions at high order, Trefftz spaces can only approximate smooth exact solutions to the governing equation at a high order;
    \item in order to achieve the same order of accuracy, Trefftz spaces require much less degrees of freedom than polynomial spaces do.
\end{itemize}
In this sense, a given Trefftz space is more efficient than a polynomial space in approximating a smaller set of functions.
%
Note that other methods, like Trefftz methods, make use of problem-dependent function spaces.
These include 
Embedded Trefftz methods \cite{LSEmb23},
Discontinuous-Petrov-Galerkin methods \cite{BOTT,DEMK},
multiscale methods \cite{Altmann_Henning_Peterseim_2021,Peters}
and
XFEM methods \cite{ABDELAZIZ20081141,STROUBOULIS20014081,MELENK1996289}.

Implementing these Trefftz methods requires bases of exact solutions to the governing PDE. For harmonic wave propagation in homogeneous media, such basis can be constructed from plane waves, from spherical waves, from evanescent waves, or from Bessel functions. However, for most PDE-driven problems, no exact solutions are known in closed form to construct local bases. In the past, this fact has limited the extent of application of Trefftz methods.

As for quasi-Trefftz DG methods, they were precisely introduced to overcome this limitation while retaining the efficiency of problem-dependent function spaces, like-Trefftz spaces.
Their specificity is to relax the property of being an exact local solution, the so-called Trefftz property, into the property of being an approximate local solution, a so-called quasi-Trefftz property. 

Quasi-Trefftz spaces have been investigated for scalar problems. To take advantage of the oscillating behavior, the first ansatz of quasi-Trefftz functions was wave-based, the idea being to add higher-degree polynomial terms in the phase of a plane wave. Other wave-based ansatz were also introduced, as well as purely a polynomial ansatz, see for instance \cite{3fam}. 
More recently in \cite{imbertgerard2025scalQTspaces}, polynomial quasi-Trefftz spaces and bases were tackled systematically for problems governed by scalar PDEs, thanks to an abstract linear algebra framework, 
for operators of arbitrary order ($\geq 1$) and in any dimension ($\geq 2$). This work relies on the fact that, in the scalar case, the principal part of the differential operator (i.e. terms of higher order derivatives with constant coefficients) is surjective between spaces of homogeneous polynomials with associated orders.
By contrast,  the system case is fundamentally different from the scalar case: for a vector-valued differential operator, the corresponding principal part is not necessarily surjective between appropriate spaces of homogeneous polynomials. Hence in this case the systematic approach developed for the scalar is not sufficient to handle even simple problems such as a first-order system for Helmholtz.

Both for time-domain and frequency-domain models, 
first-order formulations for wave propagation problems can be leveraged to develop new numerical methods efficient in the particularly challenging high-frequency regime \cite{MP21,MODAVE2023112459,PescumaPhD,RIVET2024109187,RivetPhD}.
Following this direction, this work is a stepping stone towards the development of quasi-Trefftz methods for more general problems governed by Friedrichs systems \cite{KOF,Ern2021}.
It tackles a problem governed by a first-order formulation of the Helmholtz equation.
While a quasi-Trefftz method for a first-order system formulation of the wave equation was already studied in \cite{IMS}, there the notion of quasi-Trefftz space was only defined via a quasi-Trefftz space for an auxiliary second-order scalar equation. 
With a new approach, the present work is concerned with the study of quasi-Trefftz spaces for a first-order system formulation, but instead of relying on an auxiliary scalar equation, it focuses solely on the first-order system. This is crucial for further development of quasi-Trefftz methods, 
as a decoupling approach is not always possible for more general first-order systems.

\subsection{The PDE problem and corresponding notation}\label{ssec:PDEnot}
Given a harmonic frequency $\omega$ (in $s^{-1}$), a constant density $\rho$ (in $kg\cdot m^2$), and a variable sound velocity $c$ (in $m \cdot s^{-1}$), consider the following first-order variable-coefficient partial differential system for the scalar pressure $p$ (in $N\cdot m^{-1}$) and the vector-valued velocity $\bv$ (in $m \cdot s^{-1}$):
\begin{equation}
\label{Helmholtz_fosNEW}
\left\{
\begin{array}{rcl}
\displaystyle \i \omega\rho  \bv &=& \nabla p,\\
\displaystyle \frac{\i \omega}{\rho c(x)^2} p &=& \nabla\cdot \bv.
\end{array}\right.
\end{equation}
Moreover, consider the following second-order variable-coefficient Helmholtz for the scalar pressure $p$:
\begin{equation}
\label{Helmholtz_soeNEW}
-\Delta p -\frac{\omega^2}{c(x)^2}p = 0.
\end{equation}
Then it is clear that
\begin{itemize}
    \item for any pair $(p,\bv)$ satisfying \eqref{Helmholtz_fosNEW}, then $p$ satisfies \eqref{Helmholtz_soeNEW}, and
    \item for any $p$ satisfying \eqref{Helmholtz_soeNEW}, then the pair $(p,\frac{1}{\i \omega\rho} \nabla p)$ satisfies \eqref{Helmholtz_fosNEW}.
\end{itemize}
In this sense, \eqref{Helmholtz_fosNEW} is a first-order-system formulation of the Helmholtz equation \eqref{Helmholtz_soeNEW}. Due to their variable coefficients, in general, neither the system nor the equation have exact solutions in closed form.

To allow for a more concise presentation, this work tackles the problem in dimension 2. Therefore the problem at stake has three scalar components for three scalar unknowns $(p,v_x,v_y)$:
\begin{equation}\label{eq:syst}
\left\{
\begin{array}{rcl}
\displaystyle \i \omega\rho  v_x &=& \partial_x p,\\
\displaystyle \i \omega\rho  v_y &=& \partial_y p,\\
\displaystyle \frac{\i \omega }{\rho c(x)^2} p &=& \partial_x v_x+\partial_y v_y.
\end{array}\right.
\end{equation}
Corresponding scalar differential operators are then defined for convenience:

\begin{equation*}
\begin{array}{rccl}
\mathcal L &:& p \in \mathcal C^{2} &\mapsto -\Delta p-\frac{\omega^2}{c(x)^2}p, \\
\mathcal S_1 &:& (p,v_x,v_y)\in \mathcal C^{1}\times\left(\mathcal C^{0}\right)^2 &\mapsto -\i \omega\rho  v_x + \partial_x p, \\
\mathcal S_2&:& (p,v_x,v_y)\in\mathcal C^{1}\times\left(\mathcal C^{0}\right)^2 &\mapsto -\i \omega\rho  v_y + \partial_y p, \\
\mathcal S_3 &:&(p,v_x,v_y)\in\mathcal C^{1}\times\left(\mathcal C^{0}\right)^2 &\mapsto -\frac{\i \omega}{\rho c(x)^2} p + \partial_x v_x+\partial_y v_y.
\end{array}
\end{equation*}


\subsection{Polynomial spaces and corresponding  notation}
The graded structure of polynomial spaces has been the corner stone in the study of polynomial quasi-Trefftz spaces for scalar PDEs. This will also be the case for systems of PDEs.

In this work, for any $k\in\mathbb N_0$\footnote{Throughout this work with use the notation $\mathbb N_0 = \mathbb N \cup\{0\}$.}, the space of polynomials of degree at most equal to $k$ will be denoted $\mathbb P^k$ and the space of homogeneous polynomials of degree equal to $k$ will be denoted $\widetilde{\mathbb P}^k$. Hence the graded structure can be described by
$$
\forall \degp\in\mathbb N_0,\
\mathbb P^\degp = \bigoplus_{k=0}^\degp \widetilde{\mathbb P}^k.
$$
Moreover, given any function $f\in \mathbb P^\degp$,
$f=\displaystyle\sum_{k=0}^{\degp} F^k$ with $F^k\in \widetilde{\mathbb P}^k$ for all $k$ from $0$ to $\degp$.

\subsection{Quasi-Trefftz properties and corresponding notation}
Defining the notion of quasi-Trefftz property fundamentally requires to choose a notion of approximate solution. As in previous works, the local notion of approximate solution at the root of this work is defined as a Taylor expansion approximation. In the context of a quasi-Trefftz DG method, these Taylor expansions would be defined locally on each element of the mesh. While this is later illustrated thanks to numerical results in Section \ref{sec:NR}, without loss of generality, the focus of the theoretical part of this work is limited to local notions at a given point. Given a point $\x_0\in\mathbb R^2$ and a non-negative integer $q$, the operator acting on a function of class $\mathcal C^q$ at $\x_0$ to compute its truncated Taylor polynomial of degree $q$ will be denoted hereafter $\mathcal T_q$.

In the case of a scalar equation of the form $\mathcal D u=0$, a given function $\varphi$ satisfies a quasi-Trefftz property for example if the Taylor expansion of its image by $\mathcal D$ is zero up to some order:
$$
    \mathcal T_q\mathcal D\varphi = 0.
$$
Besides, the various ansatz studied in the literature all rely on a polynomial, and corresponding quasi-Trefftz spaces are defined for a fixed maximal degree of the polynomial, denoted hereafter $\degp$, and an order of the quasi-Trefftz property that depends both on $\degp$ and on the order $\gamma$ of the differential operator, namely $\degp-\gamma$.
As a few reminders on previous work on scalar PDEs:
\begin{itemize}
    \item a phase-based GPW is a function $\varphi:\x\mapsto \exp P(\x)$ for some polynomial $P$, given an integer $d$ a corresponding quasi-Trefftz space can defined as
    $$
        \{\varphi:\x\mapsto \exp P(\x), P\in\mathbb P^\degp,\mathcal T_{\degp-\gamma}\mathcal D\varphi = 0  \},
    $$
    and this space has infinite dimension; these were first introduced in \cite{10.1093/imanum/drt030} and further studied in \cite{lmig_interpolation};
    \item an amplitude-based GPW is a function $\varphi:\x\mapsto P(\x)\exp (\i \omega \mathbf k \cdot \x)$ for some polynomial $P$, unit vector $\mathbf k$, and frequency $\omega$, given an integer $d$ and a real $\omega$ a corresponding quasi-Trefftz space can defined as
    $$
        \{\varphi:\x\mapsto P(\x)\exp (\i \omega \mathbf k \cdot \x), P\in\mathbb P^\degp, \mathbf k\in\mathbb C^2, |\mathbf k|=1,\mathcal T_{\degp-\gamma}\mathcal D\varphi = 0  \},
    $$
    and this space also has infinite dimension; these were first introduced in \cite{doi:10.1137/20M136791X} and further studied in \cite{3fam};
    \item a polynomial quasi-Trefftz function is a function $\varphi:\x\mapsto P(\x)$ for some polynomial $P$, given an integer $d$ a corresponding quasi-Trefftz space can defined as
    $$
        \{\varphi:\x\mapsto P(\x), P\in\mathbb P^\degp,\mathcal T_{\degp-\gamma}\mathcal D\varphi = 0  \},
    $$
    and this is a finite-dimensional space; 
    a recent comprehensive study of these \cite{LMIGpol}.
\end{itemize}
For the sake of clarity, this case is restricted to the choice of a polynomial ansatz. A further comment on the GPW case can be found in Section \ref{sec:GPWnHelmcomment}.

Compared to the scalar case, the case of a system of equations requires additional choices. Since both the ansatz and the system have more than one scalar component, both the polynomial degree and the notion of quasi-Trefftz property have to be defined component-wise.

\begin{remark}
    Such a quasi-Trefftz notion relying on a Taylor expansion only makes sense if the variable coefficients of the differential system (or equation) are smooth enough for their derivatives up to the desired order to exist. All results in this work will include hypothesis guaranteeing that this is the case.
\end{remark}

For the sake of compactness, the Taylor expansion centered at $\x_0$ of the differential operator's variable coefficient $\frac1{c^2}$ will be denoted
$$
    \mathcal T_q \dfrac{1}{c^2}=\kappa =\sum_{k=0}^{q} \kappa^k 
     \text{ with }
     \forall k\in[\![0,q]\!], \kappa^k \in\widetilde{\mathbb P}^k.
$$

\section{Quasi-Trefftz spaces}
In order to explore the computational aspects of the construction of quasi-Trefftz bases, this section first introduces a definition of a quasi-Trefftz space, and then proposes a different formulation of this space. This will then be the starting point to propose two explicit algorithms for the computation of quasi-Trefftz bases, presented in the next section. Finally it provides important approximation properties of this quasi-Trefftz space.

\subsection{Pressure-velocity definition}
The notion of quasi-Trefftz space is introduced here considering the first order system as a starting point. As a reminder, this is a local definition at a given point in $\mathbb R^2$.
\begin{definition}\label{def:qTspace}
    Given 
    a point $\x_0\in\mathbb R^2$,
    a polynomial degree $\degp\in\mathbb N$ with $\degp\geq 2$,
    as well as the parameters of the partial differential operators $(\mathcal S_1,\mathcal S_2,\mathcal S_3)$ defining system \eqref{eq:syst}, namely
    $\omega\in\mathbb R$,
    $\rho\in\mathbb R_+$, and 
    $c\in\mathcal C^{\degp-2}$ at $\x_0$ with $c(\x_0)\neq 0$, the corresponding quasi-Trefftz space is defined as
    $$
        \Sspace_\degp:=
        \left\{
        (p,v_x,v_y)\in\mathbb P^{\degp}\times\left(\mathbb P^{\degp-1}\right)^2 \quad\Big|\quad
    \left\{
    \begin{array}{l}
        \mathcal T_{\degp-1}\big(\mathcal S_1(p,v_x,v_y)\big)=0 \\
        \mathcal T_{\degp-1}\big(\mathcal S_2(p,v_x,v_y)\big)=0 \\
        \mathcal T_{\degp-2}\big(\mathcal S_3(p,v_x,v_y)\big)=0
    \end{array}
    \right. \qquad
    \right\}.
    $$
\end{definition}
Two important choices have been made in this definition:
\begin{itemize}
    \item 
    the polynomial degree of the three scalar unknowns are not the same; more precisely, the difference between the degree of the pressure component $p$ and the degree of the velocity components $(v_x,v_y)$ is equal to one;
    \item 
    the order of the quasi-Trefftz property for the three components of the system are not the same either, and they are chosen, following previous work on scalar problems, as the degree of the (highest-order) derivative term in the truncated quantity.
\end{itemize}
\begin{remark}\label{rmk:qvsexact}
    Since the scalar differential operators $(\mathcal S_1,\mathcal S_2)$ have constant coefficients, then 
    $$
        (p,v_x,v_y)\in\mathbb P^{\degp}\times\left(\mathbb P^{\degp-1}\right)^2 
        \Rightarrow
        \mathcal S_i(p,v_x,v_y) \in \mathbb P^{\degp-1} \ \forall i\in\{1,2\},
    $$
    therefore the first two scalar quasi-Trefftz properties in the definition of $\Sspace_\degp$ are in this case exact Trefftz properties. They can then be compactly written as $-\i \omega\rho\bv+\nabla p=0$.
\end{remark}

\subsection{Decoupled pressure}
The strong relation existing between the Helmholtz equation and the first-order system of interest will be fundamental to prove the following result, leveraging a decoupled (Helmholtz) equation for the pressure component.
\begin{proposition}\label{prop:qtspace}
    Given 
    a point $\x_0\in\mathbb R^2$,
    a polynomial degree $\degp\in\mathbb N$ with $\degp\geq 2$,
    as well as the parameters of the Helmholtz partial differential operator $\mathcal L$ introduce in Section \ref{ssec:PDEnot}, namely
    $\omega\in\mathbb R$,
    $\rho\in\mathbb R_+$, and 
    $c\in\mathcal C^{\degp-2}$ at $\x_0$ with $c(\x_0)\neq 0$, then the corresponding quasi-Trefftz space satisfies
    $$
        \Sspace_\degp=
        \left\{(p,v_x,v_y)\in\mathbb P^{\degp}\times\left(\mathbb P^{\degp-1}\right)^2 \quad|\quad \mathcal T_{\degp-2}\big(\mathcal L(p)\big)=0, \; \i \omega\rho  \bv = \nabla p\right\}.
    $$
\end{proposition}
In other words, even though the space $\Sspace_\degp$ introduced in Definition \ref{def:qTspace} is defined directly from the quasi-Trefftz approximation of the first-order system, it can also be expressed from a quasi-Trefftz approximation of a decoupled Helmholtz equation for the pressure and a post-processing to get the velocity.
\begin{proof}
($\subset$)   
Consider $(p,v_x,v_y)\in\Sspace_\degp$.
Then, according to Remark \ref{rmk:qvsexact}, the two first scalar quasi-Trefftz properties in the definition of $\Sspace_\degp$ read 
$$S_1(p,v_x,v_y)=0\text{ and }S_2(p,v_x,v_y)=0,$$
or equivalently $\i \omega\rho  \bv = \nabla p$. As a consequence, the third quasi-Trefftz property can be conveniently expressed as
$$
    \mathcal T_{\degp-2}\big(\mathcal S_3(p,v_x,v_y)\big)=0
    \Leftrightarrow
    \mathcal T_{\degp-2}\left(\frac{\i \omega}{\rho c(x)^2} p\right) =\nabla \left( \frac1{\i \omega\rho} \nabla p \right).
$$
Equivalently, since $\omega\rho$ is constant and $p\in\mathbb P^\degp$, it can also be expressed as follows
$$
    \mathcal T_{\degp-2}\left(-\frac{\omega^2}{ c(x)^2} p -\Delta p\right)=0.
$$
This proves the first inclusion by definition of the Helmholtz operator $\mathcal L$.

($\supset$) 
Consider $p\in\mathbb P^{\degp}$ satisfying the quasi-Trefftz property $\mathcal T_{\degp-2}\big(\mathcal L(p)\big)=0$ and $\bv :=\frac 1{\i \omega\rho} \nabla p$.
Following a similar argument, since $\omega\rho$ is constant and $\Delta p\in\mathbb P^{\degp-2}$, the quasi-Trefftz property can be conveniently expressed as follows
$$
    \mathcal T_{\degp-2}\big(\mathcal L(p)\big)=0
    \Leftrightarrow
    \mathcal T_{\degp-2}\left(\frac{\omega^2}{ c(x)^2} p\right) =\Delta p
    \Leftrightarrow
    \mathcal T_{\degp-2}\left(\frac{\i \omega}{\rho c(x)^2} p\right) =\nabla \bv,
$$
which in turn is equivalent to $\mathcal T_{\degp-2}\big(\mathcal S_3(p,v_x,v_y)\big)=0$.
Moreover, since $p\in\mathbb P^{\degp}$ and $\omega\rho$ is constant, then $\frac 1{\i \omega\rho} \nabla p\in\left(\mathbb P^{\degp-1}\right)^2$ and so does $\bv$ by definition.
According to Remark \ref{rmk:qvsexact}, this is equivalent to
$$\mathcal T_{\degp-1}(S_1(p,v_x,v_y))=0\text{ and }\mathcal T_{\degp-1}(S_2(p,v_x,v_y))=0.$$
This proves the second inclusion, and therefore completes the proof of the proposition.
\end{proof}

\subsection{Best approximation properties}
Proving high-order convergence of Galerkin-type methods relies on best approximation properties of discrete spaces.
In the case of quasi-Trefftz spaces, such properties can be proved via a well-chosen Taylor polynomial of the function to be approximated, as in the scalar case \cite{LMIGpol}. 

\begin{proposition}\label{prop:TuinQT}
    Given 
    a point $\x_0\in\mathbb R^2$,
    a polynomial degree $\degp\in\mathbb N$ with $\degp\geq 2$,
    as well as the parameters of the partial differential operators $(\mathcal S_1,\mathcal S_2,\mathcal S_3)$ defining system \eqref{eq:syst}, namely
    $\omega\in\mathbb R$,
    $\rho\in\mathbb R_+$, and 
    $c\in\mathcal C^{\degp-2}$ at $\x_0$ with $c(\x_0)\neq 0$,
    for any $(p,v_x,v_y)\in\mathcal C^{\degp+1}\times \mathcal C^{\degp}\times \mathcal C^{\degp}$ such that $\mathcal S_i (p,v_x,v_y)=0$ in a neighborhood of $\x_0$ for all $i\in\{1,2,3\}$, the Taylor polynomial of order $\degp$ of $(p,v_x,v_y)$ belongs to the quasi-Trefftz space of order $\degp$: $\mathcal T_{\degp}(p,v_x,v_y)\in\Sspace_\degp$. 
\end{proposition}
\begin{remark}
Of course, in general, differential operators and Taylor truncation operators do not commute, since for instance $\mathcal T_0\circ\partial_x\neq\partial_x\circ\mathcal T_0$ as the images of $f:\x=(x,y)\mapsto x$ are different ($1 \neq 0$).
Yet as we will see below, some such operators do commute when applied to smooth solution to the differential system. This relies on the fact that given $m\in\mathbb N$ and a multi-index $(\alpha_1,\alpha_2)\in\mathbb N_0^2$ with $\alpha_1+\alpha_2 \leq m-1$ then $\mathcal T_{m-\alpha_1-\alpha_2}\circ[\partial_x^{\alpha_1}\partial_y^{\alpha_2}] = [\partial_x^{\alpha_1}\partial_y^{\alpha_2}]\circ\mathcal T_m$.
\end{remark}
\begin{proof}
Since $(p,v_x,v_y)\in\mathcal C^{\degp+1}\times \mathcal C^{\degp}\times \mathcal C^{\degp}$, then 
$$
    \left\{
    \begin{array}{l}
        \mathcal S_1(p,v_x,v_y)=\mathcal S_1(\mathcal T_{\degp}p,\mathcal T_{\degp-1}v_x,\mathcal T_{\degp-1}v_y) + \mathcal S_1( p-\mathcal T_{\degp}p,v_x-\mathcal T_{\degp-1}v_x,v_y-\mathcal T_{\degp-1}v_y ),\\
        \mathcal S_2(p,v_x,v_y)=\mathcal S_2(\mathcal T_{\degp}p,\mathcal T_{\degp-1}v_x,\mathcal T_{\degp-1}v_y) + \mathcal S_2( p-\mathcal T_{\degp}p,v_x-\mathcal T_{\degp-1}v_x,v_y-\mathcal T_{\degp-1}v_y ) ,\\
        \mathcal S_3(p,v_x,v_y)=\mathcal S_3(\mathcal T_{\degp}p,\mathcal T_{\degp-1}v_x,\mathcal T_{\degp-1}v_y) + \mathcal S_3( p-\mathcal T_{\degp}p,v_x-\mathcal T_{\degp-1}v_x,v_y-\mathcal T_{\degp-1}v_y ).
    \end{array}
    \right.
$$
But
$$
    \left\{
    \begin{array}{l}
        \mathcal T_{\degp-1}\big(\mathcal S_1(\mathcal T_{\degp}p,\mathcal T_{\degp-1}v_x,\mathcal T_{\degp-1}v_y)\big) \in\mathbb P^{\degp-1}, \\
        \mathcal T_{\degp-1}\big(\mathcal S_2(\mathcal T_{\degp}p,\mathcal T_{\degp-1}v_x,\mathcal T_{\degp-1}v_y)\big) \in\mathbb P^{\degp-1}, \\
        \mathcal T_{\degp-2}\big(\mathcal S_3(\mathcal T_{\degp}p,\mathcal T_{\degp-1}v_x,\mathcal T_{\degp-1}v_y)\big) \in\mathbb P^{\degp-2},
    \end{array}
    \right. 
    \text{ while }
    \left\{
    \begin{array}{l}
        \mathcal S_1(p-\mathcal T_{\degp}p,v_x-\mathcal T_{\degp-1}v_x,v_y-\mathcal T_{\degp-1}v_y) \in \ker \mathcal T_{\degp-1} ,\\
        \mathcal S_2(p-\mathcal T_{\degp}p,v_x-\mathcal T_{\degp-1}v_x,v_y-\mathcal T_{\degp-1}v_y) \in\ker \mathcal T_{\degp-1}, \\
        \mathcal S_3(p-\mathcal T_{\degp}p,v_x-\mathcal T_{\degp-1}v_x,v_y-\mathcal T_{\degp-1}v_y) \in\ker \mathcal T_{\degp-2},
    \end{array}
    \right. 
$$
hence
$$
    \left\{
    \begin{array}{l}
        \mathcal T_{\degp-1}\mathcal S_1(p,v_x,v_y)=\mathcal T_{\degp-1}\mathcal S_1(\mathcal T_{\degp}p,\mathcal T_{\degp-1}v_x,\mathcal T_{\degp-1}v_y) ,\\
        \mathcal T_{\degp-1}\mathcal S_2(p,v_x,v_y)=\mathcal T_{\degp-1}\mathcal S_2(\mathcal T_{\degp}p,\mathcal T_{\degp-1}v_x,\mathcal T_{\degp-1}v_y)  ,\\
        \mathcal T_{\degp-2}\mathcal S_3(p,v_x,v_y)=\mathcal T_{\degp-2}\mathcal S_3(\mathcal T_{\degp}p,\mathcal T_{\degp-1}v_x,\mathcal T_{\degp-1}v_y) .
    \end{array}
    \right.
$$
Assuming that $(p,v_x,v_y)$ is an exact solution of the differential system, then all three components of the left-hand-side are zero. This concludes the proof.
\end{proof}

As a direct consequence, polynomial quasi-Trefftz spaces $\Sspace_\degp$ enjoy the following best approximation property. 
\begin{proposition}\label{prop:BAP}
    Given 
    a point $\x_0\in\mathbb R^2$,
    a polynomial degree $\degp\in\mathbb N$ with $\degp\geq 2$,
    as well as the parameters of the partial differential operators $(\mathcal S_1,\mathcal S_2,\mathcal S_3)$ defining system \eqref{eq:syst}, namely
    $\omega\in\mathbb R$,
    $\rho\in\mathbb R_+$, and 
    $c\in\mathcal C^{\degp-2}$ at $\x_0$ with $c(\x_0)\neq 0$,
    for any $(p,v_x,v_y)\in\mathcal C^{\degp+1}\times \mathcal C^{\degp}\times \mathcal C^{\degp}$ such that $\mathcal S_i (p,v_x,v_y)=0$ in a neighborhood of $\x_0$ for all $i\in\{1,2,3\}$ 
    there exists a quasi-Trefftz function $(p_a,v_{x,a},v_{y,a}) \in\Sspace_\degp$ such that, in a neighborhood of $\x_0$, 
    $$
    \left\{\begin{array}{l} 
        \|(p,v_x,v_y)(\x)- (p_a,v_{x,a},v_{y,a})(\x)\|_{\mathbb C^3}\leq C \|\x-\x_0\|_{\mathbb C^3}^{\degp+1}\\
        \|(\nabla p,\nabla v_x,\nabla v_y)(\x)- (\nabla p_a,\nabla v_{x,a},\nabla v_{y,a})(\x)\|_{\mathbb C^3}\leq C \|\x-\x_0\|_{\mathbb C^3}^{\degp}
    \end{array}\right.
$$
for some constant $C$ independent of $\mathbf x$, but depending on $(p,v_x,v_y)$, $\mathbf x_0$ and $\degp$, where $\|\cdot\|_{\mathbb C^3}$ is the euclidean norm  on ${\mathbb C^3}$.
\end{proposition}
\begin{proof}
According to Proposition \ref{prop:TuinQT}, $(\mathcal T_\degp p,\mathcal T_{\degp-1}v_{x},\mathcal T_{\degp-1}v_{y})\in\Sspace_\degp$.
The result follows from choosing $(p_a,v_{x,a},v_{y,a})=(\mathcal T_\degp p,\mathcal T_{\degp-1}v_{x},\mathcal T_{\degp-1}v_{y})$.
\end{proof}

\section{Two explicit algorithms}
Thanks to Definition \ref{def:qTspace} and Proposition \ref{prop:qtspace}, two distinct explicit algorithms can be derived for the construction of individual quasi-Trefftz functions.
In both cases, thanks to the graded structure of polynomial spaces, the construction of a (polynomial) quasi-Trefftz function can be performed one homogeneous component at a time and iteratively by increasing the homogeneous component's degree.

Any element $(p,v_x,v_y)\in\Sspace_\degp$ is comprised of three polynomials expressed in terms of their respective homogeneous components, namely $\{ P^k,k\in[\![0,\degp]\!] \}$, $\{ V_x^k,k\in[\![0,\degp-1]\!] \}$ and $\{ V_y^k,k\in[\![0,\degp-1]\!] \}$, as
$$
     p=\displaystyle\sum_{k=0}^{\degp} P^k,\
     v_x=\displaystyle\sum_{k=0}^{\degp-1} V_x^k \text{ and }
     v_y=\displaystyle\sum_{k=0}^{\degp-1} V_y^k,
$$
and $\V^k$ will denote $(V_x^k,V_y^k)$.
In turn, for any given degree $k$, these homogeneous components are expressed in terms of their respective coefficients in the canonical basis, namely $\{ \lambda^k_\ell,\ell\in[\![0,k]\!] \}$, $\{ \mu_\ell^k,\ell\in[\![0,k]\!] \}$ and $\{ \eta_\ell^k,\ell\in[\![0,k]\!] \}$, as
$$
    P^k(x,y)=\sum_{\ell=0}^k \lambda^k_\ell x^\ell y^{k-\ell},\
    V_x^k(x,y)=\sum_{\ell=0}^k \mu^k_\ell x^\ell y^{k-\ell}\text{ and }
    V_y^k(x,y)=\sum_{\ell=0}^k \eta^k_\ell\;x^\ell y^{k-\ell}.
$$

The name of the game is then to identify a sequence of well-posed linear systems that can be solved iteratively as long as all right-handside terms are computed in previous iterations.

For a given degree $\degp\in\mathbb N$ with $\degp\geq 2$, the two explicit algorithms presented in this section both construct a vector-valued function belonging to the quasi-Trefftz space $\Sspace_\degp$  given $2\degp+1$ complex values.
Fixing these values is referred to as the initialization of the algorithm.

\subsection{From the coupled definition}
Starting from the definition of a quasi-Trefftz space provided in Definition \ref{def:qTspace} and Remark \ref{rmk:qvsexact}, a convenient way to express that $(p,v_x,v_y)\in\mathbb P^{\degp}\times\left(\mathbb P^{\degp-1}\right)^2$ belongs to $\Sspace_\degp$ is 
$$
    \left\{\begin{array}{l}
        \nabla p =\i \omega\rho\bv ,\\\displaystyle
        \nabla\cdot \bv = \frac{\i \omega}{\rho }\mathcal T_{\degp-2}\left(\frac{p}{c^2} \right),
    \end{array}\right.
    \Leftrightarrow
    \left\{
    \begin{array}{l}
        \forall k\in[\![0,\degp-1]\!],\
        \nabla P^{k+1} = \i \omega\rho \V^k,\\\displaystyle
        \forall k\in[\![0,\degp-2]\!],\
        \nabla\cdot \V^{k+1} = \frac{\i \omega}{\rho }\sum_{m=0}^k\kappa^m P^{k-m},
    \end{array}
    \right.
$$
where the equivalence simply comes from the graded structure of polynomial spaces.
Constructing an element in $\Sspace_\degp$ means constructing a solution $(p,v_x,v_y)\in\mathbb P^{\degp}\times\left(\mathbb P^{\degp-1}\right)^2$ of these equivalent systems.

The second system contains for each value of $k$ from $0$ to $\degp -2$ two first-order differential equations, and an additional first-order differential equation for $k = \degp-1$.
As in the scalar case, a solution could be constructed one homogeneous component at a time for increasing values of $k$, as system can be rewritten as 
\begin{equation}\label{sys:FOS}
    \left\{
    \begin{array}{l}
    \forall k\in[\![0,\degp-2]\!],\    
    \left\{
    \begin{array}{l}
        \nabla P^{k+1} = \mathbf F^k \text{ with }\mathbf F^k:=\i \omega\rho \V^k,\\\displaystyle
        \nabla\cdot \V^{k+1} = G^k \text{ with }G^k:=\frac{\i \omega}{\rho }\sum_{m=0}^k\kappa^m P^{k-m},
    \end{array}
    \right.    \\
    \nabla P^{\degp} = \mathbf F^{\degp-1} \text{ with }\mathbf F^{\degp-1}:=\i \omega\rho \V^{\degp-1},
    \end{array}
    \right.
\end{equation}
where at each iteration $k$ the right-hand side would be known from previous iterations, as long as each equation has a solution.
Unlike in the scalar case though, one of the differential operators at stake here, namely the gradient operator acting between spaces of homogeneous polynomials, is not surjective.
Indeed, for example, there exists no polynomial $p\in\widetilde{\mathbb P}^2$  such that $\partial_x p = 2y$ and $\partial_y p = 0$.
However, the divergence operator acting between spaces of homogeneous polynomials is surjective, and provides some flexibility to solve the desired problem: it is possible to construct a solution to equations $\nabla\cdot \V^{k+1} = G^k$ such that the solution $\V^{k+1}$ belongs to the range of the gradient operator. This can be evidenced as follows.

While any vector-valued constant polynomial is in the range of the gradient operator, the situation is different for vector-valued homogeneous polynomials of higher degree.
For any $k\in[\![0,\degp-2]\!]$ in two dimensions the system's equations expressed in terms of polynomial coefficients can be conveniently combined as follows
$$
    \left\{
    {\renewcommand{\arraystretch}{1.75}%
    \begin{array}{l}\displaystyle
        \nabla\cdot \V^{k+1} = G^k \Rightarrow
        \forall \ell \in [\![ 0 , k ]\!],
        (\ell+1) \mu^{k+1}_{\ell+1} + ({k+1}-\ell) \eta^{k+1}_\ell = G^k_{\ell},
        \\\displaystyle
        \nabla P^{k+2} = \i \omega\rho \V^{k+1}\Rightarrow 
        \forall \ell \in [\![ 0 , k ]\!], 
        \dfrac{1 }{ (k+1-\ell)} \eta^{k+1}_{\ell+1} 
        =  \dfrac{1}{(\ell+1)} \mu^{k+1}_\ell,
    \end{array}
    }\right.
$$
\begin{equation}\label{eq:Vk+1}
    \Rightarrow
    \left\{
    {\renewcommand{\arraystretch}{1.75}%
    \begin{array}{l}\displaystyle
        \forall \ell \in [\![ 0 , k-1 ]\!],
        \frac{(\ell+1)(\ell+2)}{k-\ell} \eta^{k+1}_{\ell+2} + ({k+1}-\ell) \eta^{k+1}_\ell = G^k_{\ell},
        \\\displaystyle
        \mu^{k+1}_{k+1}  = \frac 1{k+1}(G^k_{k}- \eta^{k+1}_k)
        \\\displaystyle
        \forall \ell \in [\![ 0 , k ]\!], 
        \mu^{k+1}_\ell = \dfrac{(\ell+1) }{ (k+1-\ell)} \eta^{k+1}_{\ell+1} .
    \end{array}
    }\right.
\end{equation}
Hence, first, for any given values of $\left(\eta^{k+1}_0,\eta^{k+1}_1\right)$, the values $\left\{\eta^{k+1}_\ell, 2\leq\ell\leq k+1\right\}$ can be computed by the first equation; second, the values $\left\{\mu^{k+1}_\ell, 0\leq\ell\leq k+1\right\}$ can be computed by the second and third equations. 

It is then easy to justify that the $\V^{k+1}$ constructed like this belongs to the range of the gradient operator. Indeed, the gradient of any polynomial $P$ defined by
\begin{equation}\label{eq:solvegrad}
    P(x,y)=\frac{\mu^{k+1}_{k+1}}{k+2}x^{k+2}+\sum_{\ell=0}^{k+1} \frac{1}{k+2-\ell} \eta_{\ell}^{k+1} x^\ell y^{k+2-\ell} 
\end{equation}
is equal to 
$$
    \left\{\begin{array}{l}\displaystyle
        \partial_x P(x,y)=\mu^{k+1}_{k+1}x^{k+1}+\sum_{\ell=1}^{k+1} \frac{\ell}{k+2-\ell} \eta_{\ell}^{k+1} x^{\ell-1} y^{k+2-\ell}
        \\\displaystyle
        \partial_y P(x,y)= \sum_{\ell=0}^{k+1} \eta_{\ell}^{k+1} x^\ell y^{k+1-\ell} 
    \end{array}\right.
$$
But because of how the coefficients of $V_x$ are constructed from the coefficients of $V_y$ in \eqref{eq:Vk+1}, as expected, this is equivalent to stating that $\nabla P = \V^{k+1}$, and thus $\V^{k+1}$ defined in \eqref{eq:Vk+1} is in the range of the gradient. 

Turning back to the construction of a solution to System \eqref{sys:FOS}:
\begin{itemize}
    \item \eqref{eq:Vk+1} constructs a solution  to $\nabla\cdot \V^{k+1} = G^k$, and
    \item \eqref{eq:solvegrad} constructs a solution to $\nabla \left( \frac 1{\i \omega\rho } P^{k+2}\right)=\mathbf V^{k+1}$.
\end{itemize}
As for the remaining equation, that is $\nabla P^1=i\omega\rho\V^0$, it is trivial to verify that the formula corresponding to \eqref{eq:solvegrad} reads $P^1(x,y) = i\omega\rho (\mu^0_0 x+\eta^0_0 y)$.
Altogether, this can be implemented as described in Algorithm \ref{algo_fosNEW}. 
\begin{center}
\begin{algorithm}[H]
\caption{$(p,v_x,v_y) \gets$ {\tt ConstructIndivCoupled}($\degp$, $\{(\alpha_{k},\beta_{k}),0\leq k\leq \degp-1\}$, $\gamma$,$\{\kappa_k, 0\leq k\leq \degp-2 \}$,$\rho$,$\omega$)}
\label{algo_fosNEW}

\SetKwFunction{FSolGrad}{SolveGradient}
\SetKwFunction{FSolDiv}{SolveDivergence}
\SetKwProg{Fn}{Function}{:}{}

\KwIn{$\degp\in\mathbb N$ with $\degp\geq 2$, $\{(\alpha_{k},\beta_{k}),0\leq k\leq \degp-1\}\in(\mathbb C^2)^{\degp}$, $\gamma\in\mathbb C$.}
\KwOut{$p\in\mathbb P^{\degp}$,$(v_x,v_y)\in(\mathbb P^{\degp-1})^2$}

$(\lambda^0_0,\mu^0_0,\eta^0_0)\gets (\gamma,\alpha_0,\beta_0)$ \\
\For{$k = 0$ \KwTo $\degp-2$}{
$ \left\{\lambda^{k+1}_\ell, 0\leq\ell\leq k+1\right\}$
$\gets  $
\FSolGrad{$k$,$\left\{(\mu^{k}_\ell,\eta^{k}_\ell), 0\leq\ell\leq k\right\}$} \\
    \For{$\ell = 0$ \KwTo $k$}{
        $\displaystyle G^k_\ell \gets \i\omega(\sum_{m=0}^{k} \;\sum_{j=\max\{0,m+\ell-k\}}^{\min\{m,\ell\}} \kappa^{m}_j \lambda^{k-m}_{\ell-j})/\rho$ \label{line:fos_betaNEW}
    }
    $\left\{(\mu^{k+1}_\ell,\eta^{k+1}_\ell), 0\leq\ell\leq k+1\right\}\gets$\FSolDiv{$k$,$\left\{G^{k}_\ell, 0\leq\ell\leq k\right\}$, $(\alpha_{k+1},\beta_{k+1})$}
}

$ \left\{\lambda^{\degp}_\ell, 0\leq\ell\leq \degp\right\}$
$\gets  $
\FSolGrad{$\degp-1$,$\left\{(\mu^{\degp-1}_\ell,\eta^{\degp-1}_\ell), 0\leq\ell\leq \degp-1\right\},\omega,\rho$}

\Return 
$\displaystyle p(x,y) =\sum_{k=0}^{\degp}\sum_{\ell=0}^k \lambda^k_\ell x^\ell y^{k-\ell}$,
$\displaystyle v_x(x,y) =\sum_{k=0}^{\degp-1}\sum_{\ell=0}^k \mu^k_\ell x^\ell y^{k-\ell}$,
$\displaystyle v_y(x,y) =\sum_{k=0}^{\degp-1}\sum_{\ell=0}^k \eta^k_\ell x^\ell y^{k-\ell}$

    \Fn{\FSolGrad{$k$,$\left\{(\mu^{k}_\ell,\eta^{k}_\ell), 0\leq\ell\leq k\right\}\in(\mathbb C^2)^{k+1}$ coefficients of $F_x$ and $F_y$,$\omega$,$\rho$}}{

        \For{$\ell = 0$ \KwTo $k$}
        {
            $\lambda^{k+1}_\ell\gets \i\omega\rho\eta^{k}_\ell/(k+1-\ell)$ 
        }
    $\lambda^{k+1}_{k+1}\gets \i\omega\rho\mu^{k}_{k}/(k+1)$ 
    
        \Return $ \left\{\lambda^{k+1}_\ell, 0\leq\ell\leq k+1\right\}$ 
}

    \Fn{\FSolDiv{$k$,$\left\{G^{k}_\ell, 0\leq\ell\leq k\right\}\in\mathbb C^{k+1}$ coefficients of $G^k$, $(\alpha,\beta)\in\mathbb C^2$}}{

        $(\eta^{k+1}_{0},\eta^{k+1}_{1})\gets (\alpha,\beta)$

        \For{$\ell = 0$ \KwTo $k-1$}
        {
            $\eta^{k+1}_{\ell+2}\gets (k-\ell)/((\ell+1)(\ell+2))  ( G^k_\ell-(k+1-\ell)\eta^{k+1}_\ell )$ 
        }
        \For{$\ell = 0$ \KwTo $k$}
        {
            $\mu^{k+1}_{\ell}\gets (\ell+1)/(k+1-\ell)   \eta^{k+1}_{\ell+1} $ 
        }
    $\mu^{k+1}_{k+1}\gets (G^{k}_{k}-\eta^{k+1}_k)/(k+1)$ 
    
        \Return $\left\{(\mu^{k+1}_\ell,\eta^{k+1}_\ell), 0\leq\ell\leq k+1\right\}\in(\mathbb C^2)^{k+2}$ coefficients of solution
}
\end{algorithm}
\end{center}

\subsection{From the decoupled pressure}

Starting from Proposition \ref{prop:qtspace}, another convenient way to express that $(p,v_x,v_y)\in\mathbb P^{\degp}\times\left(\mathbb P^{\degp-1}\right)^2$ belongs to $\Sspace_\degp$ is
$$
    \left\{
    \begin{array}{l}\displaystyle
         \Delta p = -\omega^2\mathcal T_{\degp-2}\left( \frac{p}{c^2}\right),\\\displaystyle
        \bv = \frac{1}{\i \omega\rho}\nabla p ,
    \end{array}\right.
    \Leftrightarrow
    \left\{
    \begin{array}{l}
        \displaystyle
        \forall k\in[\![0,\degp-2]\!],\
         \Delta P^{k+2} = -\omega^2\sum_{m=0}^{\degp-2}\kappa^m P^{k-m},\\\displaystyle
        \forall k\in[\![0,\degp-1]\!],\
        V^{k} = \frac{1}{\i \omega\rho}\nabla P^{k+1} ,
    \end{array}
    \right.
$$
where the equivalence simply comes from the graded structure of polynomial spaces.
Here again, constructing an element in $\Sspace_\degp$ means constructing a solution $(p,v_x,v_y)\in\mathbb P^{\degp}\times\left(\mathbb P^{\degp-1}\right)^2$ of these equivalent systems.

The second system contains for each value of $k$ from $0$ to $\degp-2$ one second-order scalar differential equation for the pressure together with an explicit formula for the velocity in terms of pressure, and one additional explicit formula for $k=\degp-1$. Again, a solution could be constructed one homogeneous component at a time for increasing values of $k$, as the system can be rewritten as 
\begin{equation}\label{sys:SOS}
    \left\{
    \begin{array}{l}
        \displaystyle
        \forall k\in[\![0,\degp-2]\!],\  
        \Delta P^{k+2} = F^k \text{ with }F^k:=-\omega^2\sum_{m=0}^k\kappa^m P^{k-m},
        \\\displaystyle
        \forall k\in[\![0,\degp-1]\!],\  
        \V^{k} = \nabla G^{k+1} \text{ with }G^{k+1}:=\frac{1}{\i \omega\rho }P^{k+1}.
    \end{array}
    \right.
\end{equation}
Here the contruction of the pressure and velocity can be completely decoupled. The first set of scalar differential equations has a solution at each iteration $k$, it falls under the theory developed in \cite{LMIGpol} so at each iteration $k$ the equation has a solution since the right-hand side would be known from previous iterations. Then the second set of equations provides an explicit solution.
Altogether, this can be implemented as described in
Algorithm \ref{algo_sosNEW}.

\begin{algorithm}[H]
\caption{$(p,v_x,v_y) \gets$ {\tt ConstructIndivDecoupled}($\degp$, $\{(\alpha_{k},\beta_{k}),0\leq k\leq \degp-1\}$, $\gamma$,$\{\kappa_k, 0\leq k\leq \degp-2 \}$,$\rho$,$\omega$)}
\label{algo_sosNEW}

\SetKwFunction{FCompGrad}{CompGrad}
\SetKwFunction{FSolLap}{SolveLap}
\SetKwProg{Fn}{Function}{:}{}

\KwIn{$\degp\in\mathbb N$ with $\degp\geq 2$, $\{(\alpha_{k},\beta_{k}),0\leq k\leq \degp-1\}\in(\mathbb C^2)^{\degp}$, $\gamma\in\mathbb C$.}
\KwOut{$p\in\mathbb P^{\degp}$,$(v_x,v_y)\in(\mathbb P^{\degp-1})^2$}

$(\lambda^0_0,\lambda^1_0,\lambda^1_1)\gets (\gamma,\alpha_0,\beta_0)$ \\
\For{$k = 2$ \KwTo $\degp$}{
    \For{$\ell = 0$ \KwTo $k-2$}{
        $\displaystyle F^k_\ell \gets -\omega^2\sum_{m=0}^{k-2} \;\sum_{j=\max\{0,m+\ell+2-k\}}^{\min\{m,\ell\}} \kappa^{m}_j \lambda^{k-2-m}_{\ell-j}$ 
    }

    $ \left\{\lambda^{k}_\ell, 0\leq\ell\leq k\right\}$
$\gets  $
\FSolLap{$k$,$\left\{F^{k-2}_\ell, 0\leq\ell\leq k-2\right\}$,$(\alpha_{k-2},\beta_{k-2})$}
}
\For{$k = 0$ \KwTo $\degp-1$}{
    $\left\{(\mu^{k}_\ell,\eta^{k}_\ell), 0\leq\ell\leq k\right\}\gets $\FCompGrad{$k$,$\left\{\lambda^{k+1}_\ell/(\i\omega\rho), 0\leq\ell\leq k+1\right\}$}
}

\Return 
$\displaystyle p(x,y) =\sum_{k=0}^{\degp}\sum_{\ell=0}^k \lambda^k_\ell x^\ell y^{k-\ell}$,
$\displaystyle v_x(x,y) =\sum_{k=0}^{\degp-1}\sum_{\ell=0}^k \mu^k_\ell x^\ell y^{k-\ell}$,
$\displaystyle v_y(x,y) =\sum_{k=0}^{\degp-1}\sum_{\ell=0}^k \eta^k_\ell x^\ell y^{k-\ell}$

    \Fn{\FSolLap{$k$,$\left\{F^{k-2}_\ell, 0\leq\ell\leq k-2\right\}\in\mathbb C^{k-1}$ coefficients of $F^{k-2}$, $(\alpha,\beta)\in\mathbb C^2$}}{

        $(\lambda^{k}_{0},\lambda^{k}_{1})\gets (\alpha,\beta)$

        \For{$\ell = 0$ \KwTo $k-2$}
        {
            $\lambda^k_{\ell+2}\gets\big(F^{k-2}_\ell-(k-\ell)(k-\ell-1)\;\lambda^k_\ell\big)/\big((\ell+2)(\ell+1)\big)$ 
        }
    
        \Return $\left\{(\lambda^{k}_\ell), 0\leq\ell\leq k\right\}\in(\mathbb C^2)^{k+1}$ coefficients of solution
}

   \Fn{\FCompGrad{$k$,$\left\{G^{k+1}_\ell, 0\leq\ell\leq k+1\right\}\in\mathbb C^{k+2}$ coefficients of $G^{k+1}$}}{
        \For{$\ell = 0$ \KwTo $k$}{
            $\mu^k_\ell \gets (\ell+1) G^{k+1}_{\ell+1}$ 
            \\
            $\eta^k_\ell \gets (k+1-\ell) G^{k+1}_{\ell}$
        }
        \Return $\left\{(\mu^{k}_\ell,\eta^{k}_\ell), 0\leq\ell\leq k\right\}$
   }

\end{algorithm}

\subsection{Extension to 3D}
While the present work focuses on System \eqref{Helmholtz_fosNEW} in two-dimensional setting, it is a natural question to ask how would it extend to three dimensions. Thanks to the general structure of Algorithms \ref{algo_fosNEW} and \ref{algo_sosNEW}, it is actually straightforward to provide an answer. 

First, in terms of notation, we now have $\bv=(v_x,v_y,v_z)$.
Any element $(p,\bv)\in \mathbb P^{\degp}\times\left(\mathbb P^{\degp}\right)^3$ is comprised of four polynomials expressed in terms of their respective homogeneous components, namely $\{ P^k,k\in[\![0,\degp]\!] \}$, $\{ V_x^k,k\in[\![0,\degp-1]\!] \}$, $\{ V_y^k,k\in[\![0,\degp-1]\!] \}$, and $\{ V_z^k,k\in[\![0,\degp-1]\!] \}$, as
$$
     p=\displaystyle\sum_{k=0}^{\degp} P^k,\
     v_x=\displaystyle\sum_{k=0}^{\degp-1} V_x^k,
     v_y=\displaystyle\sum_{k=0}^{\degp-1} V_y^k \text{ and }
     v_z=\displaystyle\sum_{k=0}^{\degp-1} V_z^k.
$$
In turn, for any given degree $k$, each homogeneous component is expressed in terms of its coefficients in the canonical basis, namely $\{ \lambda^k_{j,\ell},j,\ell\in[\![0,k]\!], 0\leq j+\ell\leq k \}$, $\{ \mu_{j,\ell}^k,j,\ell\in[\![0,k]\!], 0\leq j+\ell\leq k \}$, $\{ \eta_{j,\ell}^k,j,\ell\in[\![0,k]\!], 0\leq j+\ell\leq k\}$ and $\{ \nu_{j,\ell}^k,j,\ell\in[\![0,k]\!], 0\leq j+\ell\leq k\}$, as
$$
    P^k(x,y,z)=\sum_{0\leq j+\ell\leq k} \lambda^k_{j,\ell} x^j y^\ell z^{k-j-\ell},\
    V_x^k(x,y,z)=\sum_{0\leq j+\ell\leq k} \mu^k_{j,\ell} x^j y^\ell z^{k-j-\ell},
$$
$$
    V_y^k(x,y,z)=\sum_{0\leq j+\ell\leq k} \eta^k_{j,\ell} x^j y^\ell z^{k-j-\ell} \text{ and }
    V_z^k(x,y,z)=\sum_{0\leq j+\ell\leq k} \nu^k_{j,\ell} x^j y^\ell z^{k-j-\ell}.
$$
Moreover, $\V^k$ will denote $(V_x^k,V_y^k,V_z^k)$

To generalize Algorithm \ref{algo_fosNEW}, first notice that System \eqref{sys:FOS} is valid as is in three dimensions. The aim is then to construct a solution to the equation $\nabla \cdot \V^{k+1} = G^k$ such that $\V^{k+1}$ lies in the range of the gradient operator.
Given that 
$$
    \left\{
    {\renewcommand{\arraystretch}{1.75}%
    \begin{array}{ll}\displaystyle
        \nabla\cdot \V^{k+1} = G^k \Rightarrow&
        \forall 0\leq j+\ell \leq k,
        (j+1) \mu^{k+1}_{j+1,\ell} + (\ell+1) \eta^{k+1}_{j,\ell+1}+ (k+1-j-\ell)\nu^{k+1}_{j,\ell} = G^k_{j,\ell},
        \\\displaystyle
        \nabla P^{k+2} = \i \omega\rho \V^{k+1}\Rightarrow& 
        \mu^{k+1}_{j,\ell} 
        =\dfrac{j+1}{(k+1-j-\ell)} \nu^{k+1}_{j+1,\ell}, 
        \eta^{k+1}_{j,\ell}=\dfrac{\ell+1}{(k+1-j-\ell)} \nu^{k+1}_{j,\ell+1},
    \end{array}
    }\right.
$$
imples that
\begin{equation}
\label{sys_3dNEW}
    \left\{
    {\renewcommand{\arraystretch}{1.75}%
    \begin{array}{l}
        \displaystyle \forall 0\leq j+\ell \leq k-1,  \frac{(j+1)(j+2)}{k-j-\ell} \nu^{k+1}_{j+2,\ell}+\frac{(\ell+1)(\ell+2)}{k-j-\ell} \nu^{k+1}_{j,\ell+2} + ({k+1}-j-\ell) \nu^{k+1}_{j,\ell} = G^k_{j,\ell},
        \\\displaystyle \forall 0\leq j+\ell \leq k, \mu^{k+1}_{j,\ell}  =\dfrac{j+1}{(k+1-j-\ell)} \nu^{k+1}_{j+1,\ell},\;  \eta^{k+1}_{j,\ell}=\dfrac{\ell+1}{(k+1-j-\ell)} \nu^{k+1}_{j,\ell+1},
        \\\displaystyle \forall j\in[\![0,k]\!], \mu^{k+1}_{j+1,k-j} =\dfrac{1}{(j+1)} \left( -(k-j+1) \eta^{k+1}_{j,k-j+1} - \nu^{k+1}_{j,k-j} + G^k_{j,k-j}\right),
        \\\displaystyle
        \forall 0\leq j+\ell \leq k, \eta^{k+1}_{j+1,\ell} = \dfrac{(\ell+1)}{(j+1)} \mu^{k+1}_{j,\ell+1}, 
    \end{array}
    }\right.
\end{equation}
this can be performed 
\begin{enumerate}
    \item fix any given values of $\{\nu^{k+1}_{j,\ell}, j\in\{0,1\} , 0\leq \ell \leq k+1-j\}$, 
    \\compute the values $\{\nu^{k+1}_{j,\ell}, 2\leq j \leq k+1, 0\leq \ell \leq k+1-j\}$;
    \item compute the values $\{(\mu^{k+1}_{j,\ell},\eta^{k+1}_{j,\ell}), 0\leq j \leq k+1, 0\leq \ell \leq k+1-j \}$;
    \item for any given $(\eta^{k+1}_{0,k+1},\eta^{k+1}_{1,k})$, 
    \\compute the values of $\{\mu^{k+1}_{j,k+1-j}, 0\leq j\leq k+1\}$ and  $\{\eta^{k+1}_{j,k+1-j}, 2\leq j \leq k+1\}$.
\end{enumerate}
This simply corresponds to implementing three-dimensional versions of the  {\tt SolveDivergence} and {\tt SolveGradient} functions in Algorithm \ref{algo_fosNEW}.
\begin{remark}
    In this 3D case, the resulting dimension of the quasi-Trefftz space, given by the number of free parameters in the algorithm, is  $(\degp+1)^2$.
\end{remark}

As for Algorithm \ref{algo_sosNEW}, it simply requires three-dimensional versions of the  {\tt SolveLap} and {\tt CompGrad} functions. The latter is straightforward while the former can be found in \cite{LMIGpol}.

\section{Quasi-Trefftz bases}
This section focuses on procedures to construct whole bases of quasi-Trefftz spaces. Two main types of procedures will be proposed, either computing sequentially individual basis elements via the explicit algorithms from the previous section, or computing simultaneously all basis elements as the kernel of a unique linear operator.

\subsection{Constructing individual basis elements}
Given a fixed degree $d\in\mathbb N$ with $\degp\geq 2$, an individual element of $\Sspace_\degp$ can by constructed, via either Algorithm \ref{algo_fosNEW} or Algorithm \ref{algo_sosNEW}, via a set of values for the initialization.
While distinct sets of values define distinct elements of $\Sspace_\degp$, the challenge to build a basis is to define linearly independent elements in $\Sspace_\degp$.

However, in both algorithms, the given sets of values are assigned to distinct polynomial coefficients of the constructed vector-valued quasi-Trefftz function $(p,v_x,v_y)$ expressed in the canonical monomial polynomial basis, more precisely these coefficients are:
$$
\text{either }
\{\lambda_0,\mu_0,\eta_0\}\cup \{(\eta^k_0,\eta^k_1), 1\leq k\leq \degp-1\} 
\text{ or }
\{\lambda_0\}\cup\{\lambda^k_0,\lambda^k_1, 1\leq k\leq \degp\}.
$$
This has two important consequences. On the one hand, if $2\degp+1$ sets of values form a basis of $\mathbb C^{2\degp+1}$, the corresponding $2\degp+1$ quasi-Trefftz functions with these initialization sets of values are linearly independent. On the other hand, $2\degp+1$ is the dimension of the quasi-Trefftz space $\Sspace_\degp$.

For convenience, the canonical basis of $\mathbb C^{2\degp+1}$ will be denoted $\{\mathbf e_k, 1\leq k\leq 2\degp+1\}$, that is
$$
    \forall (k,\ell) \in [\![1,2\degp+1]\!]^2,
    (\mathbf e_k)_\ell = \delta_k^\ell.
$$

As a result, a basis of the quasi-Trefftz space $\Sspace_\degp$ can be constructed by iteratively applying either Algorithm \ref{algo_fosNEW} or Algorithm \ref{algo_sosNEW} for convenient initialization values to construct individual basis elements. This is summarized in Algorithms \ref{algo_fos_basis} and \ref{algo_sos_basis}.
\begin{center}
\begin{algorithm}[H]
\caption{{\tt ConstructBasisCoupled}($\degp$,$\{\kappa_k, 0\leq k\leq \degp-2 \}$,$\rho$,$\omega$)}
\label{algo_fos_basis}

\SetKwFunction{Indiv}{ConstructIndivCoupled}

\For{$j = 1$ \KwTo $2\degp+1$}{
    $(p^j,v_x^j,v_y^j) \gets$\Indiv{$\degp$,$\{((\mathbf e_j)_{2k},(\mathbf e_j)_{2k+1}),0\leq k\leq \degp-1\}$,$(\mathbf e_j)_{2\degp+1}$,$\{\kappa_k, 0\leq k\leq \degp-2 \}$,$\rho$,$\omega$}
}

\Return $\{(p^j,v_x^j,v_y^j),1\leq j\leq 2\degp+1\}$

\end{algorithm}
\end{center}


\begin{center}
\begin{algorithm}[H]
\caption{{\tt ConstructBasisDecoupled}($\degp$,$\{\kappa_k, 0\leq k\leq \degp-2 \}$,$\rho$,$\omega$)}
\label{algo_sos_basis}

\SetKwFunction{Indiv}{ConstructIndivDecoupled}

\For{$j = 1$ \KwTo $2\degp+1$}{
    $(p^j,v_x^j,v_y^j) \gets$\Indiv{$\degp$,$\{((\mathbf e_j)_{2k},(\mathbf e_j)_{2k+1}),0\leq k\leq \degp-1\}$,$(\mathbf e_j)_{2\degp+1}$,$\{\kappa_k, 0\leq k\leq \degp-2 \}$,$\rho$,$\omega$}
}

\Return $\{(p^j,v_x^j,v_y^j),1\leq j\leq 2\degp+1\}$

\end{algorithm}
\end{center}

\subsection{Computational complexity}
\label{part_complexity_expl}
Since the two previous algorithms construct a quasi-Trefftz basis functions by calling $2\degp+1$ times an algorithm constructing an individual basis function, it is then natural to focus on the computational cost of these algorithms, namely Algorithms \ref{algo_fosNEW} and \ref{algo_sosNEW}.

The focus is on counting the number of floating point operations, or flops, namely here additions and multiplications.
No distinction is made between real and complex arithmetic.
In both of these, there are various multiplicative factors that are independent of the values of the input data, and hence can be precomputed once and for all. Therefore their computation will not be taken into account in the study of the computational complexity of each algorithm.

\begin{figure}[H]
    \centering
    \subfigure[$M=4$]{\includegraphics[width=0.45\linewidth]{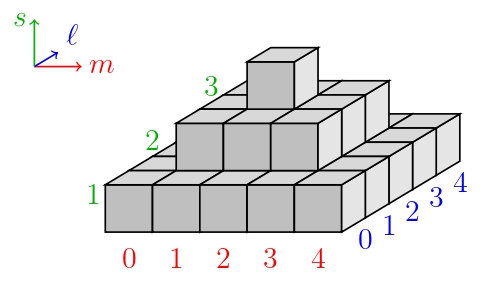}} 
    \subfigure[$M=5$]{\includegraphics[width=0.45\linewidth]{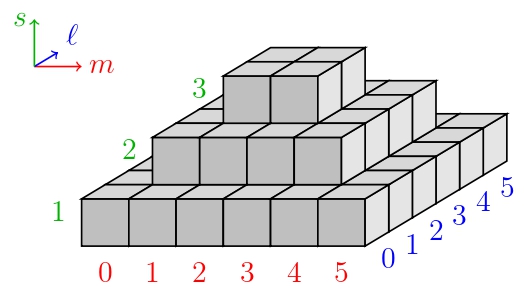}} 
    \caption{Values of the number $s=\min\{m,\ell\}-\max\{0,m+\ell-M\}+1$ of terms in the inside sum in \eqref{eq:dsum}, as a function of $(m,\ell)$ for an odd value and an even value of $M$. 
    }
    \label{fig:pyrNEW}
\end{figure}

In both algorithms the computational work is dominated by the innermost loop term
\begin{equation}\label{eq:dsum}
    \forall \ell\in[\![0,M]\!],\
    \sum_{m=0}^{M} \;\sum_{j=\max\{0,m+\ell-M\}}^{\min\{m,\ell\}} \kappa^{m}_j \lambda^{M-m}_{\ell-j}
\end{equation}
for each value of $M\in[\![ 0, \degp-2]\!]$.
For fixed values of $M$ and $\ell$, the number of terms in the sum over $j$, here denoted $s$ for convenience, depends on the values of $m$ with respect to $\ell$ and $M-\ell$ as described in Table \ref{tab:valsS}.
\begin{table}
    \begin{center}
    \begin{tabular}{|c||c|c|}
        \hline
                     & $m \leq M-\ell$ & $m > M-\ell$    \\\hline\hline
        $m\leq \ell$ & $s=m+1$         & $s=M+1-\ell$    \\\hline 
        $m> \ell$    & $s=\ell+1$      & $s=M+1-m$       \\\hline
    \end{tabular}
    \end{center}
    \caption{Values of $s=\min\{m,\ell\}-\max\{0,m+\ell-M\}+1$ depending on the signs of $m-\ell$ and $m-M+\ell$, for given values of $M$ and $\ell$.}
    \label{tab:valsS}
\end{table}
In order to find the total number of terms in \eqref{eq:dsum} for a given value of $M$, the value of $s$ can then conveniently be represented geometrically as a function of $\ell\in[\![ 0, M]\!]$ and $m\in[\![ 0, M]\!]$: one unit cube representing one term, for each pair $(m,\ell)$ the number of terms is represented by a column of cubes of height $s$.
 The volume of the corresponding pyramid, denoted $V_M$, is then precisely the total number of terms in \eqref{eq:dsum}. 
 It can then be calculated as the sum of volumes of each horizontal layer of the pyramid: 
$$
    V_M = \sum_{n=0}^{\lfloor M/2\rfloor} \left( 2n+1+ r(M) \right)^2
    \Rightarrow 
    V_{M} = \frac{(M+1)(M+2)(M+3)}{6} 
$$
where the remainder $r$ is defined as $r(M)=2(\frac M2-\lfloor \frac M2\rfloor)$ is simply zero for even values of $M$ and one for odd values.
Figure \ref{fig:pyrNEW} illustrates this geometric point of view and the difference between odd values and even values of $M$.

Moreover, in both algorithms, \eqref{eq:dsum} is computed for each value of $M$ from $0$ to $\degp-2$. Hence the total number $T_\degp$ of terms in
\begin{equation}\label{eq:dsumtot}
    \forall M\in[\![0,\degp-2]\!],\
    \forall \ell\in[\![0,M]\!],\
    \sum_{m=0}^{M} \;\sum_{j=\max\{0,m+\ell-M\}}^{\min\{m,\ell\}} \kappa^{m}_j \lambda^{M-m}_{\ell-j}
\end{equation}
can be computed as
$$
    T_\degp
    = \sum_{M=0}^{\degp-2} V_M 
    = \sum_{M=0}^{\degp-2}\frac{(M+1)(M+2)(M+3)}{6}  
     = \frac{(\degp-1)\degp(\degp+1)(\degp+2)}{24}.
$$

The total number $\tilde T_\degp$ of additions required to compute \eqref{eq:dsumtot} can be computed similarly. For $M=0$ there is a single term so no addition. 
For fixed values of $M>1$ and $\ell\in[\![0,M]\!]$, the number of additions to compute the sum over $j$ is $s-1$. Using again a geometric argument, the corresponding pyramids would be obtained from those of Figure \ref{fig:pyrNEW} by simply removing the lowest level ($s=1$), and hence their volume would be
$$
    \tilde V_M = V_M-(M+1)^2
    \Rightarrow 
    \tilde V_{M} = \frac{(M-1)M(M+1)}{6} ,
$$
and so the total number of additions would be
$$
    \tilde T_\degp
    =\sum_{M=1}^{\degp-2}\tilde V_{M} 
    =\frac {(\degp-3)(\degp-2)(\degp-1)\degp}{24} 
$$

In conclusion, calculating  \eqref{eq:dsumtot} requires a total of $T_\degp+\tilde T_\degp = (\degp-1)\degp(\degp^2-\degp+4)/12$ additions and multiplications.

The remaining computational work performed in the algorithms can be described as follows:
\begin{itemize}
    \item solving a divergence system {\tt SolveDivergence}, requiring $k+1$ additions and $3k+2$ multiplications
        performed for each $k$ from $0$ to $\degp-2$,
        so $(\degp-1)(2\degp-1)$ operations in total;
    \item solving a gradient system via {\tt SolveGradient}, requiring $k+2$ multiplications
        performed for each $k$ from $0$ to $\degp-1$,
        so $(\degp+1)(\degp+2)/2 -1$ operations in total;
    \item solving a Laplacian system via {\tt SolveLaplacian}, requiring $k-1$ additions and $2(k-1)$ multiplications
        performed for each $k$ from $2$ to $\degp$,
        so $3\degp(\degp-1)/2$ operations in total;
    \item computing a gradient via {\tt CompGrad}, requiring $2(k+1)$ multiplications
        performed for each $k$ from $0$ to $\degp-1$,
        so $\degp(\degp+1)$ operations in total.
\end{itemize}

To summarize,
the number of operations performed in Algorithm \ref{algo_fosNEW} is
$$
    \frac{(\degp-1)\degp(\degp^2-\degp+4)}{12} + (\degp-1)(2\degp-1) + \frac{(\degp+1)(\degp+2)}{2} -1
    =
    \frac{\degp^4 - 2 \degp^3 + 35 \degp^2 - 22 \degp + 12}{12}
$$
while in Algorithm \ref{algo_sosNEW} it is
$$
    \frac{(\degp-1)\degp (\degp^2-\degp+4)}{12} + \frac{3\degp(\degp-1)}2 + \degp(\degp+1)
    =\frac{\degp^4 - 2 \degp^3 + 35 \degp^2 - 10\degp}{12}
$$

To construct a quasi-Trefftz basis,
the computational complexity of Algorithm \ref{algo_fos_basis} is simply $2\degp+1$ times that of Algorithm \ref{algo_fosNEW}, and similarly
the computational complexity of Algorithm \ref{algo_sos_basis} is simply $2\degp+1$ times that of Algorithm \ref{algo_sosNEW}.
As a result, the computational complexity for the construction of a quasi-Trefftz basis from explicit algorithms can be summarized as follows.
\begin{center}
    {\renewcommand{\arraystretch}{2.5}
    \begin{tabular}{|c|c|}
        \hline
        Algorithm            & Number of {\it flops}\\\hline\hline
        \ref{algo_fos_basis} & $\displaystyle \frac{(2d+1)(\degp^4 - 2 \degp^3 + 35 \degp^2 - 22 \degp + 12)}{12}$ \\\hline
        \ref{algo_sos_basis} & $\displaystyle \frac{(2\degp+1)(\degp^4 - 2 \degp^3 + 35 \degp^2 - 10\degp)}{12}$ \\\hline
    \end{tabular}
    }
\end{center}

\subsection{Constructing all elements simultaneously}
As opposed to the sequential computation of basis functions previously discussed, it is possible to compute all the elements of the basis simultaneously.
Indeed, the quasi-Trefftz space $\Sspace_\degp$ can be expressed as the kernel of an appropriate linear operator, and a basis for the operator's kernel can be computed thanks to a Singular Value Decomposition (SVD) of any matrix representing the operator. More precisely, as a reminder, for any matrix $M\in\mathbb C^{m\times n}$ of rank $r$ and its SVD $M = U \Sigma V^*$ where the diagonal entries of $\Sigma$, $\{\sigma_k,1\leq k\leq \min(m,n)\}$ are a decreasing function of $k$, then the last $n-r$ columns of $V$ form a basis of $\ker M$. See for instance \cite{trefethen}.

In order to proceed, given
    a point $\x_0\in\mathbb R^2$,
    a polynomial degree $\degp\in\mathbb N$ with $\degp\geq 2$,
    as well as the parameters of the partial differential operator $\mathcal L$, namely
    $\omega\in\mathbb R$,
    $\rho\in\mathbb R_+$, and 
    $c\in\mathcal C^{\degp-2}$ at $\x_0$ with $c(\x_0)\neq 0$,
    define the linear operator
\begin{equation*}
    \begin{array}{rccc}
    \mathcal{Q}^F_\degp: &\mathbb P^{\degp}\times\left(\mathbb P^{\degp-1}\right)^2 &\longrightarrow &\mathbb P^{\degp-2}\times\left(\mathbb P^{\degp-1}\right)^2 \\
    &\begin{bmatrix} p \\ v_x \\ v_y \end{bmatrix} &\mapsto &\begin{bmatrix} \mathcal T_{\degp-2}(-\frac{\i \omega}{\rho c^2} p + \partial_x v_x+\partial_y v_y) \\ -\i\omega\rho v_x +\partial_x p \\ -\i\omega\rho v_y +\partial_y p \end{bmatrix}.
    \end{array}
\end{equation*}
Then the quasi-Trefftz space can indeed be expressed as $\Sspace_\degp = \Ker(\mathcal{Q}^F_\degp)$.
A matrix representing $\mathcal{Q}^F_\degp$ is uniquely defined by two bases of the domain and codomain with a given numbering of their respective elements.

\begin{figure}[H]
    \centering
    \includegraphics[width=0.5\linewidth]{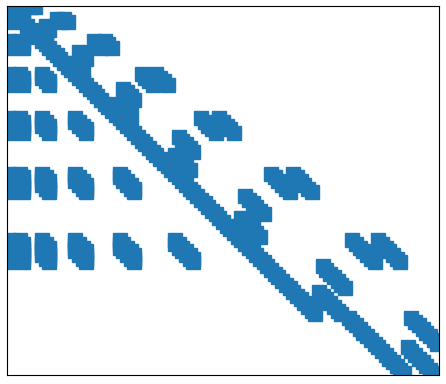}
    \caption{Structure of a matrix representing the linear operator $Q^F_\degp$ for $\degp=8$.}
    \label{mat_fNEW}
\end{figure}

For the sake of completeness, a matrix representing $\mathcal{Q}^F_\degp$ can be constructed as follows. First, both the domain and the codomain are vector spaces of vector-valued polynomials of the form
\begin{equation*}
    \mathbb P^{n}\times\left(\mathbb P^{m}\right)^2 = 
    \left\{
    \begin{array}{ll}
    \displaystyle
        \left( \bigoplus_{k=0}^{n} \mathbb H^k\times\left(\mathbb H^{k}\right)^2\right) \bigoplus 
        \left( \bigoplus_{k=n+1}^{m} \{\mathbf 0_{\mathbb P}\}\times\left(\mathbb H^{k}\right)^2\right)& n<m ,\\
    \displaystyle
        \left( \bigoplus_{k=0}^{m} \mathbb H^k\times\left(\mathbb H^{k}\right)^2\right) \bigoplus 
        \left( \bigoplus_{k=m+1}^{n} \mathbb H^{k}\times\left(\{\mathbf 0_{\mathbb P}\}\right)^2\right)& n>m,
    \end{array}
    \right.
\end{equation*}
emphasizing their graded structure. 
Defining component-wise scalar canonical polynomial basis functions as
$$
    \forall k\in\mathbb N_0,\
    \forall \ell\in[\![0,k]\!],\
    \phi^k_{\ell,1}:(x,y)\mapsto (x^\ell y^{k-\ell},0,0), 
$$
$$
    \forall k\in\mathbb N_0,\
    \forall \ell\in[\![0,k]\!],\
    \phi^k_{\ell,2}:(x,y)\mapsto (0,x^\ell y^{k-\ell},0), 
$$
$$
    \forall k\in\mathbb N_0,\
    \forall \ell\in[\![0,k]\!],\
    \phi^k_{\ell,3}:(x,y)\mapsto (0,0,x^\ell y^{k-\ell}),
$$
then consider
\begin{itemize}
    \item $\{ \phi^k_{\ell,1},k\in[\![0,\degp]\!],\ell\in[\![0,k]\!]\}\cup \{ \phi^k_{\ell,i},k\in[\![0,\degp-1]\!],\ell\in[\![0,k]\!],i\in\{2,3\} \}$, basis of $\mathbb P^{\degp}\times\left(\mathbb P^{\degp-1}\right)^2$, and
    \item $\{ \phi^k_{\ell,1},k\in[\![0,\degp-2]\!],\ell\in[\![0,k]\!]\}\cup \{ \phi^k_{\ell,i},k\in[\![0,\degp-1]\!],\ell\in[\![0,k]\!],i\in\{2,3\} \}$, basis of $\mathbb P^{\degp-2}\times\left(\mathbb P^{\degp-1}\right)^2$.
\end{itemize}
The image of each element of the domain's basis can then be expressed in the codomain's basis as follows:
\begin{itemize}
    \item 
    $\displaystyle \mathcal{Q}^F_\degp (\phi^k_{\ell,1}) = 
    -\frac{\i\omega}\rho \mathbf 1_{k<\degp-1}\sum_{m=0}^{\degp-2-k}\sum_{j=0}^{m}\kappa^m_j \phi^{m+k}_{j+\ell,1}
    +\ell \mathbf 1_{k\ell>0} \phi^{k-1}_{\ell-1,2}
    +(k-\ell)\mathbf 1_{k(k-\ell)>0} \phi^{k-1}_{\ell,3}
    $,
    \item 
    $\displaystyle \mathcal{Q}^F_\degp (\phi^k_{\ell,2}) = 
    -i\omega\rho \phi^k_{\ell,2}
    +\ell \mathbf 1_{k\ell>0}\phi^{k-1}_{\ell-1,1}
    $,
    \item 
    $\displaystyle \mathcal{Q}^F_\degp (\phi^k_{\ell,3}) =
    -i\omega\rho \phi^k_{\ell,3} 
    +(k-\ell) \mathbf 1_{k(k-\ell)>0}\phi^{k-1}_{\ell,1}
    $,
\end{itemize}
where for all $k\in\mathbb N_0$ $\{\kappa^k_j,j\in[\![0,k]\!]\}$ are the coefficients of the homogeneous polynomial $\kappa^k\in\widetilde{\mathbb P}^k$ in the canonical basis $\{\phi^k_{\ell}:(x,y)\mapsto x^\ell y^{k-\ell},\ell\in[\![0,k]\!]\}$.
In order to illustrate the sparsity pattern of the linear operator, Figure \ref{mat_fNEW} displays an example for $\degp=8$ using the following numbering of basis elements. 
For both bases, 
elements $\phi^k_{\ell,i}$ are gathered first by common values of $k$ and the corresponding sets are ordered by increasing values of $k$;
then within a constant $k$ block, they are gathered by common values of $i$ and the corresponding sets are ordered by increasing values of $i$;
finally, within a constant $k$ and $i$ block, they are ordered by increasing values of $\ell$.

\begin{figure}[H]
    \centering
    \includegraphics[width=0.5\linewidth]{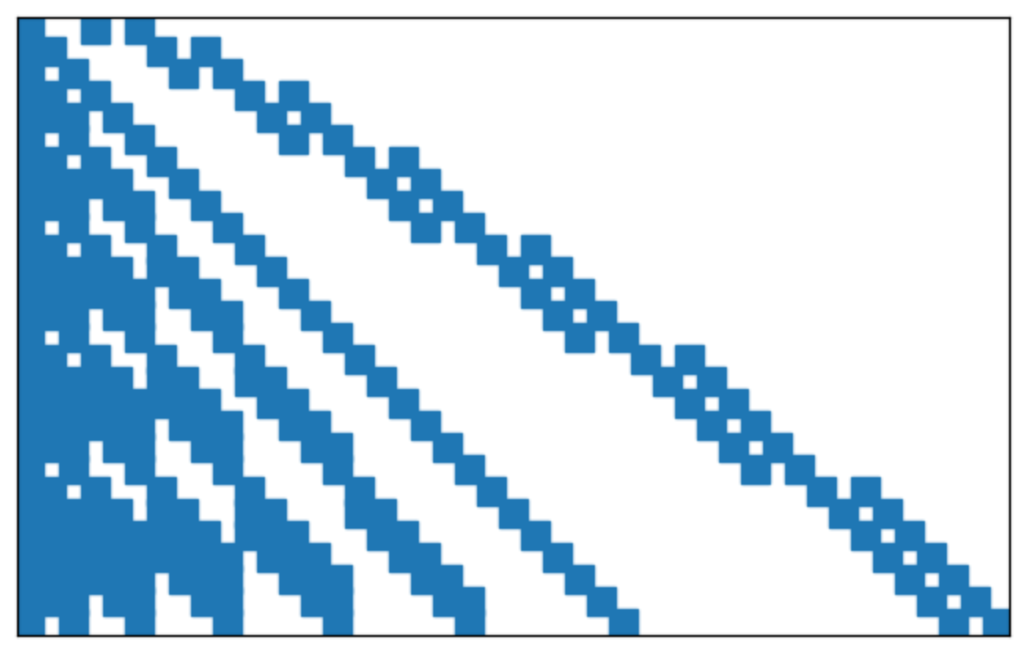}
    \caption{Structure of matrices representing the linear operators $Q^S_\degp$ and $G_\degp$ for $\degp=8$. }
    \label{mat_sNEW}
\end{figure}

An alternative, similarly relying on an SVD as opposed to an explicit algorithm, can also be derived from Proposition \ref{prop:qtspace}. 
Here, in order to proceed, given
    a point $\x_0\in\mathbb R^2$,
    a polynomial degree $\degp\in\mathbb N$ with $\degp\geq 2$,
    as well as the parameters of the partial differential operator $\mathcal L$, namely
    $\omega\in\mathbb R$,
    $\rho\in\mathbb R_+$, and 
    $c\in\mathcal C^{\degp-1}$ at $\x_0$,
    define the linear operators
\begin{equation*}
    \begin{array}{rccc}
    \mathcal{Q}^S_\degp: &\mathbb P^{\degp} &\longrightarrow &\mathbb P^{\degp-2} \\
    &p &\mapsto &\displaystyle\mathcal T_{\degp-2}\left(\Delta p + \frac{\omega^2}{c^2}p\right)
    \end{array},
    \text{ and }
    \begin{array}{rccc}
    G_\degp: &\mathbb P^{\degp} &\longrightarrow &(\mathbb P^{\degp-1})^2 \\
    &p &\mapsto &\displaystyle\frac 1{\i\omega\rho} \nabla p
    \end{array}.
\end{equation*}
The quasi-Trefftz spaces can be expressed as follows
\begin{equation*}
    \Sspace_\degp= 
    \left\{(p,v_x,v_y)\in\mathbb P^{\degp}\times\left(\mathbb P^{\degp-1}\right)^2 \quad|\quad p\in\Ker(\mathcal{Q}^S), \bv = G_\degp (p)\right\}.
\end{equation*}
Therefore the construction of a basis can be performed in two steps: first the simultaneous construction of the $p$ components as the operator's kernel can be computed thanks to an SVD of any matrix representing the $\mathcal{Q}^S_\degp$ operator, then the simultaneous construction of the $\bv$ components as the image of the $p$ components by the $G_\degp$ operator. Here the domains and codomains for both operators are scalar polynomial spaces, $\mathbb P^n$ for some value of $n$, of which the canonical bases are $\{\phi^k_{\ell}:(x,y)\mapsto x^\ell y^{k-\ell},\ell\in[\![0,k]\!],k\in[\![0,n]\!]\}$.
The image of each element of the domain's basis can then be expressed in the comdomain's basis as follows for any $(k,\ell)\in[\![0,k]\!]\times[\![0,\degp]\!]$:
$$
    \mathcal{Q}^S_\degp ( \phi^k_{\ell} )
    =
    \mathbf 1_{k>1}\mathbf 1_{  \ell>1}\ell(\ell-1)      \phi^{k-2}_{\ell-2}
    + 
    \mathbf 1_{k>1}\mathbf 1_{k-\ell>1}(k-\ell)(k-\ell-1)\phi^{k-2}_{\ell}
    + \omega^2 
    \sum_{m=0}^{\degp-2-k}\sum_{j=0}^{m}\kappa^m_j \phi^{m+k}_{j+\ell}
$$
where again for all $k\in\mathbb N_0$ $\{\kappa^k_j,j\in[\![0,k]\!]\}$ are the coefficients of the homogeneous polynomial $\kappa^k\in\widetilde{\mathbb P}^k$ in the canonical scalar polynomial basis, and
$$
    G_\degp (\phi^k_{\ell}) = \frac 1{\i\omega\rho} \left(\ell\mathbf 1_{  \ell>0}\phi^{k-1}_{\ell-1},(k-\ell)1_{k-\ell>0}\phi^{k-1}_{\ell} \right).
$$
In order to illustrate the sparsity pattern of both linear operators, Figure \ref{mat_sNEW} displays an example for $\degp=8$ using the following numbering of basis elements.
For both bases, elements are gathered by common value of $k$ and the corresponding sets are ordered by increasing values of $k$.


\begin{remark}
    Given a set of indices denoted $I$ and two vector space $\mathbb A$ and $\mathbb B$, a standard result of linear algebra states that for any $\{(x_i,y_i),i\in I \}\subset\mathbb A\oplus\mathbb B$, if $\{x_i,i\in I \}$ is a linearly independent subset of $\mathbb A$, then$\{x_i+y_i,i\in I \}$ is a linearly independent subset of $\mathbb A\oplus\mathbb B$.
    Hence the set of functions in $\mathbb P^{\degp}\times\left(\mathbb P^{\degp-1}\right)^2$ constructed above is indeed linearly independent since $\mathbb P^{\degp}\times\left(\mathbb P^{\degp-1}\right)^2 = \mathbb P^{\degp}\times\left\{(0,0)\right\}\oplus \left\{0\right\}\times\left(\mathbb P^{\degp-1}\right)^2$ and it is known from \cite{LMIGpol} that the pressure components are linearly independent.
\end{remark}

\subsection{Computational complexity}\label{ssec:compcompALGE}

The matrix $\mathcal Q^F_\degp$ has
    $\degp (3\degp+1)/2$ rows and
    $(\degp+1)(3\degp+2)/2$ columns.
Therefore, the computational cost to compute the $V$ in the SVD $\mathcal Q^F_\degp=U\Sigma V^*$  is 
$$
    2\left(\frac{\degp (3\degp+1)}2+\frac{(\degp+1)(3\degp+2)}2\right)\left(\frac{\degp (3\degp+1)}2\right)^2
    =\frac{\degp^2(27\degp^4+45 \degp^3+30\degp^2+9\degp+1)}2.
$$
A specialized SVD algorithm taking advantage of the sparsity pattern of the matrix could also be used, and this could be computationally more efficient.



As for the matrix $\mathcal Q_\degp^S$, it has
$(\degp-1)\degp/2$ rows
and 
$(\degp+1)(\degp+2)/2$ columns.
Therefore, the computational cost to compute the $V$ in the SVD $\mathcal Q^S_\degp=U\Sigma V^*$  is 
$$
    2\left( \frac{(\degp-1)\degp}2 + \frac{(\degp+1)(\degp+2)}2 \right) \left( \frac{(\degp-1)\degp}2 \right)^2
    =
    \frac{d^2(\degp^4-\degp^3-\degp+1)}2.
$$
Moreover, any matrix $Q^G_\degp$ representing the linear operator $G_\degp$ has
$
2\degp(\degp+1)/2=
\degp(\degp+1)
$ rows
and
$(\degp+1)(\degp+2)/2$ columns.
Therefore, multiplying it by a matrix with
$(\degp+1)(\degp+2)/2$ rows
and 
$2\degp+1$ columns has a computational cost of 
$$
    \left( 2\frac{(\degp+1)(\degp+2)}2-1\right)\degp(\degp+1)(2\degp+1)
    =
    \left( (\degp+1)(\degp+2)-1\right)\degp(\degp+1)(2\degp+1).
$$

\subsection{A brief comparison}
\begin{figure}[H]
    \centering
    \includegraphics[width=0.6\linewidth]{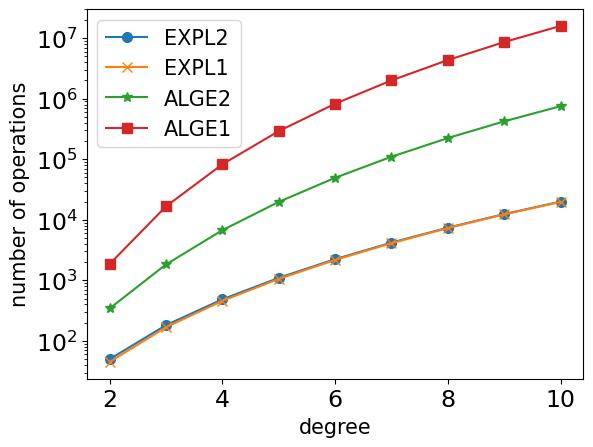}
    \caption{Number of floating-point operations required to construct one quasi-Trefftz basis of $\Sspace_\degp$, for degrees $\degp$ from 2 to 10, using the four methods (EXPL1, EXPL2, ALGE1, and ALGE2).}
    \label{fig:flops}
\end{figure}

We have just introduced four methods to build a quasi-Trefftz basis:
\begin{itemize}
    \item the explicit method based on the coupled pressure and velocity, 
    hereafter referred to as EXPL1;
    \item the explicit method based on the decoupled pressure, 
    hereafter referred to as EXPL2;
    \item the SVD method based on the coupled pressure and velocity, 
    hereafter referred to as ALGE1;
    \item the SVD method based on the decoupled pressure, 
    hereafter referred to as ALGE2.
\end{itemize}
For each method, the number of floating-point operations required is evaluated (see Subsection~\ref{part_complexity_expl} for the explicit cases and Subsection \ref{ssec:compcompALGE} for the algebraic case). Figure~\ref{fig:flops} reports these complexities for the construction of one basis of $\Sspace_\degp$, with $\degp$ ranging from 2 to 10.
The results show that both explicit methods incur very similar costs. The algebraic method associated with the decoupled problem (ALGE2) is more expensive by one to two orders of magnitude, while the algebraic method associated with the coupled problem (ALGE1) is yet another one to two orders higher.

\section{Numerical results}
\label{sec:NR}

The goal of this section is to illustrate the compared behavior of the four approaches implemented for constructing quasi-Trefftz bases, namely: EXPL1, EXPL2, ALGE1 and ALGE2,
A first simple test to verify the implementation of these four approaches is to verify that the basis functions constructed indeed satisfy the desired quasi-Trefftz properties.
A second test is to verify the best approximation properties of the space generated by the constructed basis functions.
Finally, computational costs of the different approaches can be compared.

All computations were performed in double precision on a laptop with an Intel Core i7-13700H and 32 GB RAM.

\subsection{Test cases}

The notion of quasi-Trefftz space studied in the present work 
is introduced given a polynomial degree $\degp$. Since quasi-Trefftz DG methods are high-order methods, the tests proposed hereafter systematically compare results obtained for values of $\degp$ ranging  from $2$  to $10$.

This notion 
is local, defined given a point $\x_0$.
In the context of quasi-Trefftz DG methods for PDE problems, given a meshed computational domain, quasi-Trefftz spaces are constructed locally on each element of the mesh. This will hence be done here as well: defining a test case will include defining a domain and describing a mesh of the domain, and the errors reported here will always be the worst error occuring across quasi-Trefftz spaces on each individual mesh element.

This notion 
is also introduced given three scalar differential operators $(\mathcal S_1,\mathcal S_2,\mathcal S_3)$, themselves depending on two scalar parameters $(\omega,\rho)$ and one variable coefficient $c$. Therefore defining a test case requires the value of the two parameters as well as the definition of the variable coefficient. Moreover, note that implementing the construction of bases requires the coefficients of the Taylor expansion of the variable coefficient up to order $\degp-2$ at the center of each element of the mesh.

\begin{figure}[H]
\centering
\subfigure[Test Case 1: square domain, structured triangular mesh with 800 elements, polynomial coefficient.]{\includegraphics[width=0.38\linewidth]{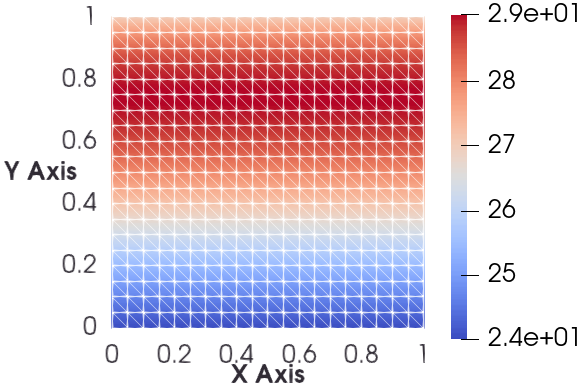}
\label{fig:test_caseNEW_1}}
\hfill
\subfigure[Test Case 2: elliptical domain around a simplified wing profile, uniform triangular mesh with 5847 elements, variable coefficient modeling temperature-induced sound speed variations.]{\includegraphics[width=0.59\linewidth]{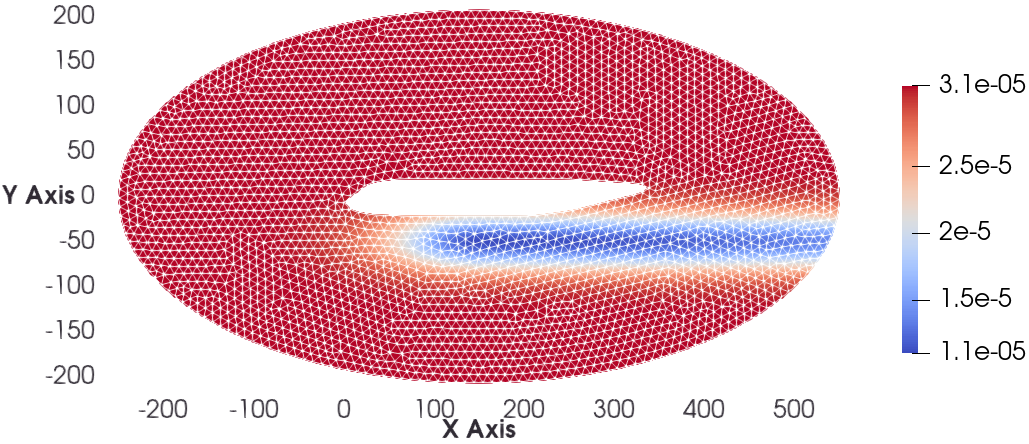}
\label{fig:test_caseNEW_2}}
\caption{Representation of the computational domains, their meshes, and the corresponding variable coefficients $\omega^2/c^2$ for the two test cases.}
\end{figure}

The test case 1 is defined as follows.
The domain is the square $\Omega=[0,1]\times[0,1]$.
It is meshed It is meshed with a structured triangular mesh consisting of 800 elements.
The variable coefficient is the polynomial $\dfrac{\omega^2}{c^2(x,y)}=24-9y^4+6y^2+6y$, see Figure \ref{fig:test_caseNEW_1}, and $(\omega,\rho)=(\pi/2,4/\pi^2)$.

The test case 2 mimics a simplified aeroacoustic application, it is defined as follows.
The domain is an ellipse surrounding a basic plane wing profile.
It is meshed with a uniform triangular mesh consisting of 5847 elements.
The variable coefficient is a toy model for variations in the sound speed induced by a heated jet. The sound speed is given by $c=\sqrt{\gamma R T}$ where $\gamma=1.4$ is the ratio of specific heats, $R\approx 287 J\cdot kg^{-1}\cdot K^{-1}$ is the specific gas constant for air, and $T$ is the temperature varying from 200 to 550 $K$. Finally, with $(\omega,\rho)=(\pi/2,4/\pi^2)$ we obtain
$$
    \dfrac{\omega^2}{c^2(x,y)}=
    \left\{\begin{array}{cc}
         3.1\times10^{-5}-2\times10^{-5}\text{e}^{-10^{-4}(x-150)^2-0.001(y+50)^2}
         & \text{ if }x\leq150, \\
         3.1\times10^{-5}-2\times10^{-5}\text{e}^{-10^{-6}(x-150)^2-0.001(y+50)^2}
         & \text{ otherwise},
    \end{array}\right.
$$
 see Figure \ref{fig:test_caseNEW_2}, and $(\omega,\rho)=(\pi/2,4/\pi^2)$.

As a reminder, given a test case, bases of local quasi-Trefftz spaces are constructed on each mesh element. 

\subsection{Computational performance}

To complement the operation count analysis (see Figure~\ref{fig:flops}), we report here the actual computation times required to construct a single quasi-Trefftz basis, averaged over all elements, for the four methods and both test cases. These timings provide a practical illustration of the computational cost of each method. 

It should however be emphasized that the reported times are only indicative. They strongly depend on the present implementation and are not the result of a performance-optimized code. Consequently, they should not be interpreted as absolute benchmarks, but rather as order-of-magnitude estimates confirming the tendencies predicted by the operation counts.

\begin{figure}[H]
    \centering
    \includegraphics[width=0.6\linewidth]{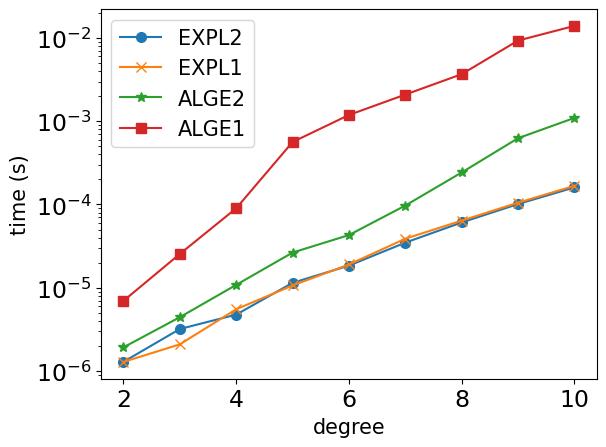}
    \caption{Computation time to construct a single quasi-Trefftz basis, averaged over all elements, comparing the four methods (EXPL1, EXPL2, ALGE2, ALGE1) for Test Case 1.}
    \label{fig:time_case1}
\end{figure}

\begin{figure}[H]
    \centering
    \includegraphics[width=0.6\linewidth]{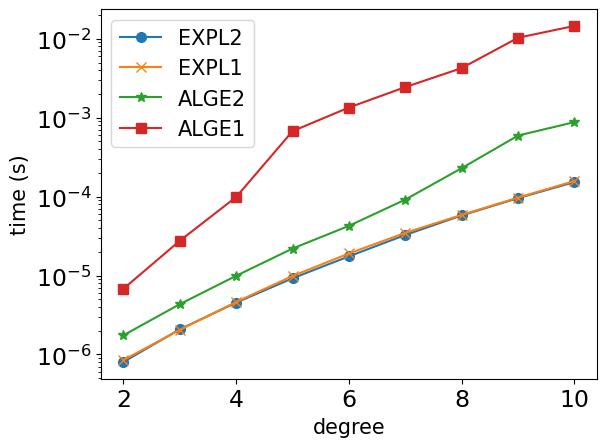}
    \caption{Computation time to construct a single quasi-Trefftz basis, averaged over all elements, comparing the four methods (EXPL1, EXPL2, ALGE2, ALGE1) for Test Case 2.}
    \label{fig:time_case2}
\end{figure}

Figures~\ref{fig:time_case1} and~\ref{fig:time_case2} summarize these timings.
As expected from the theoretical operation counts, the algebraic SVD-based approaches are more expensive than the explicit ones, with the coupled variant (ALGE1) being the most demanding, while the two explicit methods (EXPL1 and EXPL2) exhibit very similar timings. Furthermore, the same behavior is obtained for both test cases as predicted by the theory.

\subsection{Quasi-Trefftz properties}
The mere definition of a polynomial quasi-Trefftz space $\Sspace_\degp$ consists of three scalar quasi-Trefftz properties, namely
$$
    \mathcal T_{\degp-1}(-\i \omega\rho v_x + \partial_x p )= 0, 
    \quad 
    \mathcal T_{\degp-1}(-\i \omega\rho v_y + \partial_y p) = 0
    \quad \text{ and } \quad 
    \mathcal T_{\degp-2}\left(-\frac{\i \omega}{\rho c^2} p + \partial_x v_x+\partial_y v_y\right) = 0.
$$
Moreover, given that $(p,v_x,v_y)\in\mathbb P^{\degp}\times\left(\mathbb P^{\degp-1}\right)^2$, then the first two properties are actually Trefftz properties:
$$
    -\i \omega\rho v_x + \partial_x p = 0
    \quad \text{ and } \quad 
    -\i \omega\rho v_y + \partial_y p = 0.
$$
In other words:
\begin{itemize}
    \item the $\mathbb P^{\degp-1}$ polynomials $-\i \omega\rho v_x + \partial_x p$ and $-\i \omega\rho v_y + \partial_y p$ are both zero;
    \item the function $-\frac{\i \omega}{\rho ^2} p + \partial_x v_x+\partial_y v_y$ behaves like $|\x-\x_0|^{\degp-1}$ as $\x\to\x_0$.
\end{itemize}
In order to verify that any function indeed satisfies these properties
\begin{itemize}
    \item the $\ell^\infty$ norm of the coefficients of the $\mathbb P^{\degp-1}$ polynomials $-\i \omega\rho v_x + \partial_x p$ and $-\i \omega\rho v_y + \partial_y p$ expressed in the canonical monomial basis should be zero;
    \item the $L^\infty(|\x-\x_0|<h)$ norm of function $-\frac{\i \omega}{\rho c(x)^2} p + \partial_x v_x+\partial_y v_y$ should behave like $h^{\degp-1}$ as $h\to 0$.
\end{itemize}
\begin{table}[H]
\centering     {\renewcommand{\arraystretch}{1.2}
\begin{tabular}{|c|c|c|c|c|c|}
\hline
 & \degp & EXPL2 & EXPL1 & ALGE2 & ALGE1 \\
\hline
\hline
Test case 1 & 2 & 0 & $7.1 \; 10^{-15}$ & 0 &  $1.1 \; 10^{-15}$  \\ 
\hline
Test case 1 & 3 & 0 & $7.1 \; 10^{-15}$ & 0 &  $6.6 \; 10^{-15}$  \\ 
\hline
Test case 1 & 4 & 0 & $2.8 \; 10^{-14}$ & 0 &  $2.8 \; 10^{-15}$  \\ 
\hline
Test case 1 & 5 & 0 & $5.6 \; 10^{-14}$ & 0 &  $1.1 \; 10^{-14}$  \\ 
\hline
Test case 1 & 6 & 0 & $1.1 \; 10^{-13}$ & 0 &  $6.6 \; 10^{-15}$  \\ 
\hline
Test case 1 & 7 & 0 & $1.1 \; 10^{-13}$ & 0 &  $8.2 \; 10^{-15}$  \\ 
\hline
Test case 1 & 8 & 0 & $1.1 \; 10^{-13}$ & 0 &  $5.7 \; 10^{-14}$  \\ 
\hline
Test case 1 & 9 & 0 & $3.4 \; 10^{-13}$ & 0 &  $1.2 \; 10^{-14}$  \\ 
\hline
Test case 1 & 10 & 0 & $4.5 \; 10^{-13}$ & 0 &  $1.8 \; 10^{-14}$  \\ 
\hline
\hline
Test case 2 & 2 & 0 &  $1.4 \; 10^{-20}$ & 0 & $1.3\; 10^{-15}$  \\  
\hline
Test case 2 & 3 & 0 &  $1.4 \; 10^{-20}$ & 0 & $4.4\; 10^{-15}$  \\  
\hline
Test case 2 & 4 & 0 &  $1.4 \; 10^{-20}$ & 0 & $4.4\; 10^{-15}$  \\  
\hline
Test case 2 & 5 & 0 &  $3.6\; 10^{-15}$ & 0 & $8.7\; 10^{-15}$  \\  
\hline
Test case 2 & 6 & 0 &  $3.6\; 10^{-15}$ & 0 & $8.4\; 10^{-15}$  \\  
\hline
Test case 2 & 7 & 0 &  $2.4 \; 10^{-14}$ & 0 & $9.9 \; 10^{-15}$  \\  
\hline
Test case 2 & 8 & 0 &  $5.7 \; 10^{-14}$ & 0 & $7.7 \; 10^{-15}$  \\  
\hline
Test case 2 & 9 & 0 &  $1.1 \; 10^{-13}$ & 0 & $1.1 \; 10^{-14}$  \\  
\hline
Test case 2 & 10 & 0 &  $2.2 \; 10^{-13}$ & 0 & $1.2 \; 10^{-14}$  \\  
\hline
\end{tabular}}
\caption{Maximum $\ell^\infty$ norm of the polynomial coefficients of $-\i\omega \rho v_x + \partial_x p$ and $-\i\omega \rho v_y + \partial_y p$ (the first two quasi-Trefftz properties), taken over all basis functions and all elements, 
for polynomial degrees ranging from $2$ to $10$, 
for both test cases, and for the four basis construction methods (EXPL2, EXPL1, ALGE2, ALGE1).}
\label{tab:qtprop1}
\end{table}
Table~\ref{tab:qtprop1} reports the maximum $\ell^\infty$ norm of the polynomial coefficients associated with the first two quasi-Trefftz properties, over all basis functions and elements, for polynomial degrees from 2 to 10. As expected, the decoupled methods (EXPL2 and ALGE2) achieve exact satisfaction of these properties, yielding zero coefficients, because in these methods the velocity components $v_x$ and $v_y$ are explicitly defined as $v_x = \frac{1}{\i\omega\rho} \partial_x p$ and $v_y = \frac{1}{\i\omega\rho} \partial_y p$. In contrast, the coupled methods (EXPL1 and ALGE1) produce errors close to machine precision, confirming that these properties are accurately preserved.
We observe that as the polynomial degree increases, the error tends to grow, most likely due to the larger number of computations involved.

\begin{figure}[H]
\centering
\subfigure[EXPL2]{\includegraphics[width=0.49\linewidth]{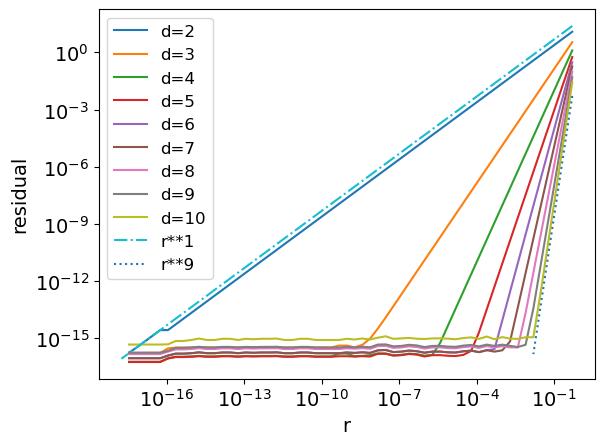}}
\subfigure[ALGE2]{\includegraphics[width=0.49\linewidth]{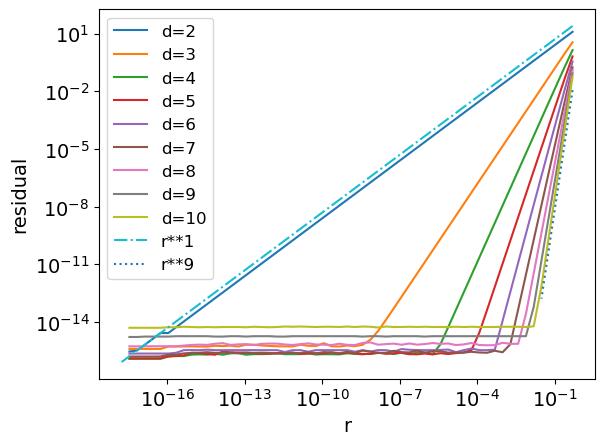}}
\subfigure[EXPL1]{\includegraphics[width=0.49\linewidth]{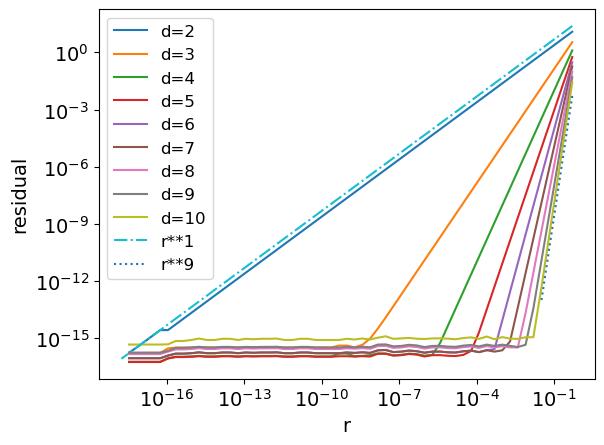}}
\subfigure[ALGE1]{\includegraphics[width=0.49\linewidth]{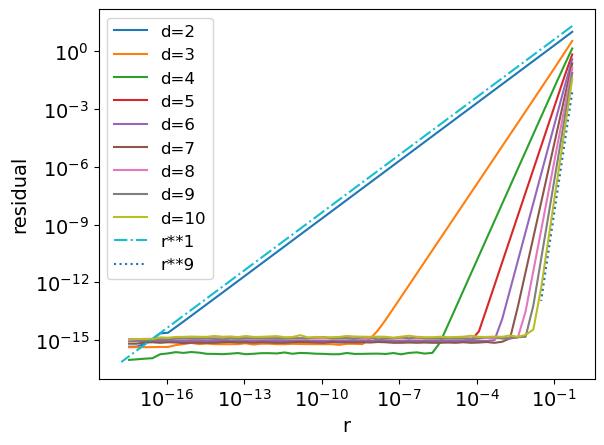}}
\caption{Maximum relative error for the residual of the third quasi-Trefftz property 
$-\frac{\mathrm{i} \omega}{\rho c(x)^2} p + \partial_x v_x + \partial_y v_y$ 
over all basis functions and elements, evaluated for decreasing distances $r$ from the element center. 
Computed for the four construction methods (EXPL2, EXPL1, ALGE2, ALGE1), for Test Case~1.}
\label{fig:qtprop1}
\end{figure}

\begin{figure}[H]
\centering
\subfigure[EXPL2]{\includegraphics[width=0.49\linewidth]{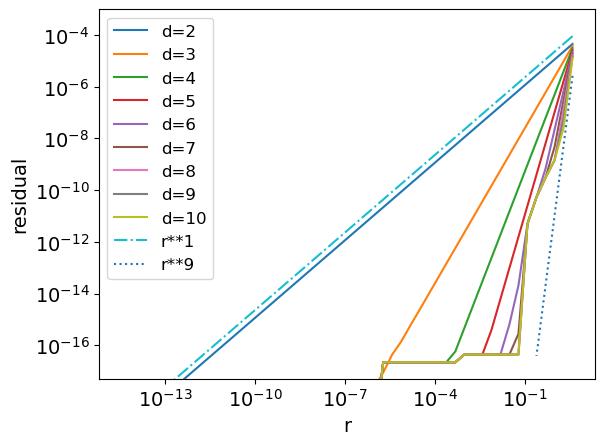}}
\subfigure[ALGE2]{\includegraphics[width=0.49\linewidth]{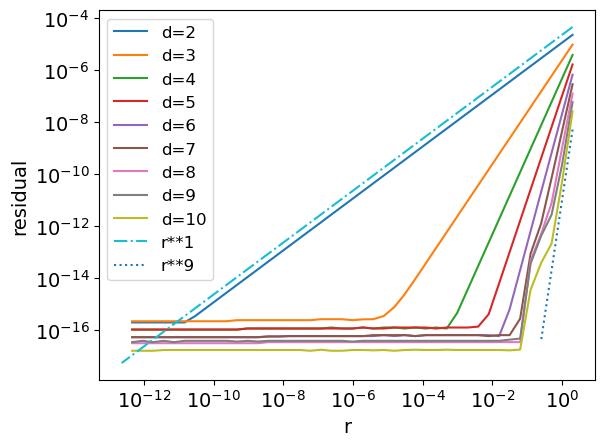}}
\subfigure[EXPL1]{\includegraphics[width=0.49\linewidth]{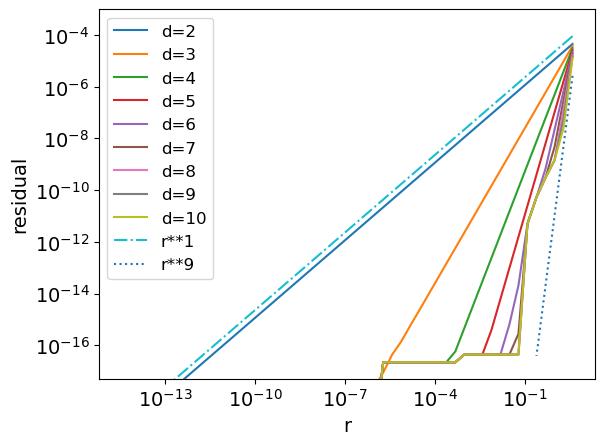}}
\subfigure[ALGE1]{\includegraphics[width=0.49\linewidth]{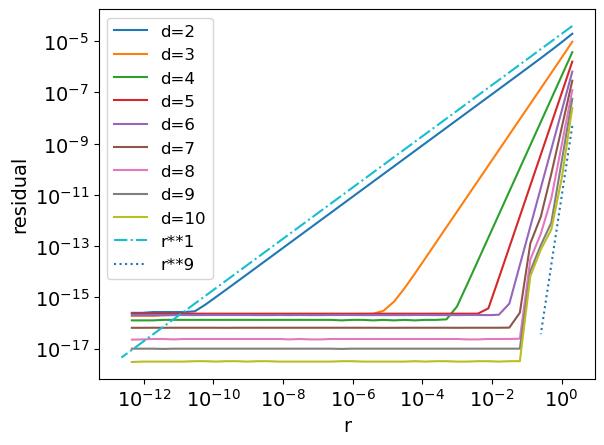}}
\caption{Maximum relative error for the residual of the third quasi-Trefftz property 
$-\frac{\mathrm{i} \omega}{\rho c(x)^2} p + \partial_x v_x + \partial_y v_y$ 
over all basis functions and elements, evaluated for decreasing distances $r$ from the element center. 
Computed for the four construction methods (EXPL2, EXPL1, ALGE2, ALGE1), for Test Case~2.}
\label{fig:qtprop2}
\end{figure}

To assess the third quasi-Trefftz property, on each element, we consider individual basis functions and for each of them we compute the maximum residual of $-\frac{\i \omega}{\rho ^2} p + \partial_x v_x+\partial_y v_y$
over points at decreasing distances $r$ from the element center. 
We estimate the $L^\infty$ norm of this residual, which according to $\Sspace_\degp$'s defintion is expected to decrease as $r^{\degp-1}$ as $r\to 0$.
We record the largest relative error among all basis functions on the element and among all elements. This procedure is performed for polynomial quasi-Trefftz spaces of degrees $\degp$ from 2 to 10. The results are shown in Figures~\ref{fig:qtprop1} and~\ref{fig:qtprop2} for the two test cases.
For Test Case 1, the expected convergence rate is clearly observed: for each polynomial degree $\degp$, the residual decays with slope $\degp-1$ until it reaches a plateau near the level of machine precision. 
This plateau appears at slightly higher values as the polynomial degree increases, which is consistent with the fact that there are more computations involved as the degree increases. 
Concerning Test Case 2, for low polynomial degrees, the results are essentially the same as in Test Case 1, with a decay of slope $\degp-1$ followed by a plateau near machine precision.
Here, this plateau appears at slightly lower values as the polynomial degree increases, which is the opposite of what was observed in the previous case. 
Moreover, for higher degrees -starting from $\degp=6$ with the explicit methods and $\degp=7$ with the algebraic ones- we observe the appearance of bumps in the error curves. When comparing the error on all elements of the mesh, these bumps only occur at the discontinuity of the derivatives of the Helmholtz coefficient $\dfrac{1}{c^2}$ while on all other elements the behavior is very similar to that observed in Test Case 1. 
This is consistent with the hypothesis of the theoretical results.


\subsection{Best approximation properties}
An exact solution to the differential system for test case 1 is given by:
\begin{equation*}
    p_{ex}(x,y)=e^{5\i x -y^3+y}, \quad  \text{ and } \qquad \bv_{ex}=\dfrac{1}{\i\omega\rho}\nabla p_{ex}.
\end{equation*}

\begin{figure}[H]
\centering
\subfigure[EXPL2]{\includegraphics[width=0.49\linewidth]{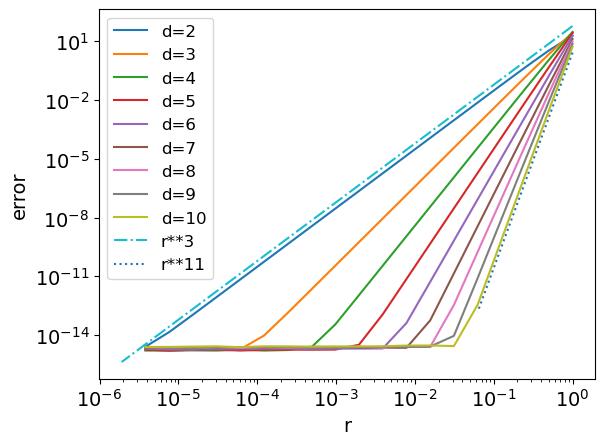}}
\subfigure[ALGE2]{\includegraphics[width=0.49\linewidth]{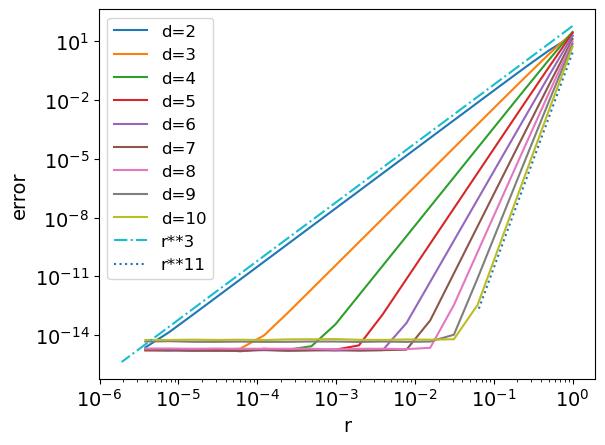}}
\subfigure[EXPL1]{\includegraphics[width=0.49\linewidth]{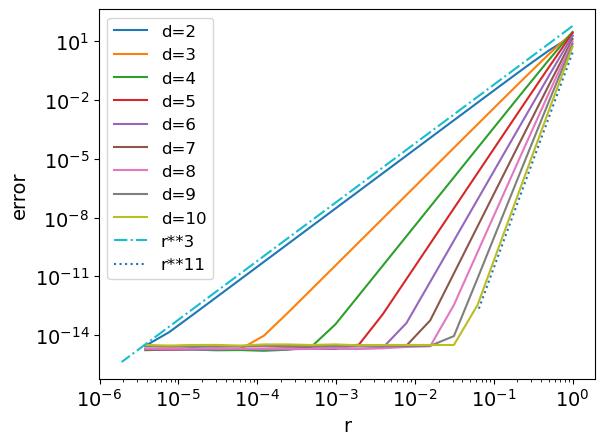}}
\subfigure[ALGE1]{\includegraphics[width=0.49\linewidth]{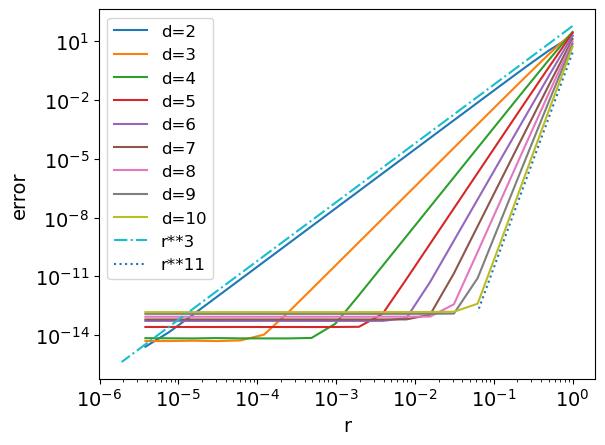}}
\caption{Error between the exact solution and quasi-Trefftz approximations for the pressure component, as a function of the distance from the element center, Test Case 1.}
\label{fig:approxp}
\end{figure}

\begin{figure}[H]
\centering
\subfigure[EXPL2]{\includegraphics[width=0.49\linewidth]{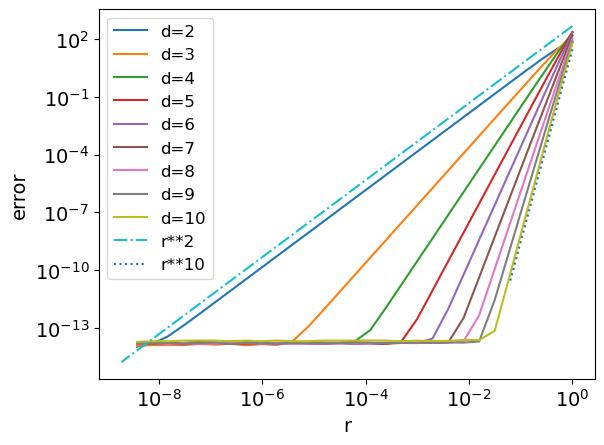}}
\subfigure[ALGE2]{\includegraphics[width=0.49\linewidth]{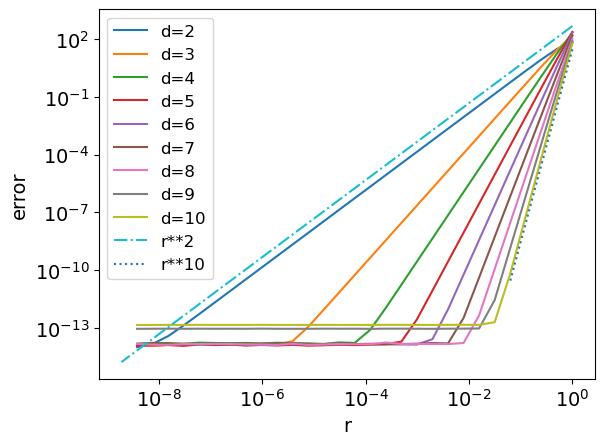}}
\subfigure[EXPL1]{\includegraphics[width=0.49\linewidth]{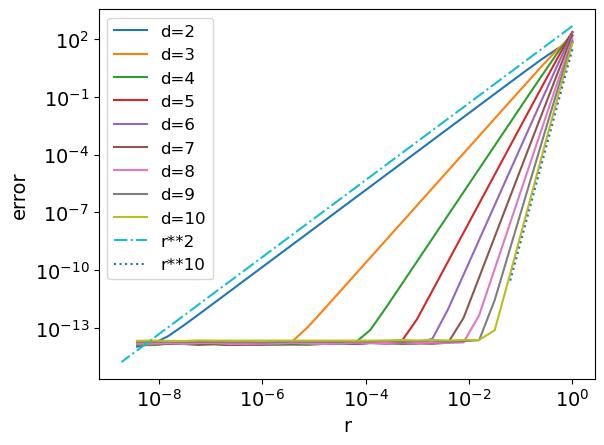}}
\subfigure[ALGE1]{\includegraphics[width=0.49\linewidth]{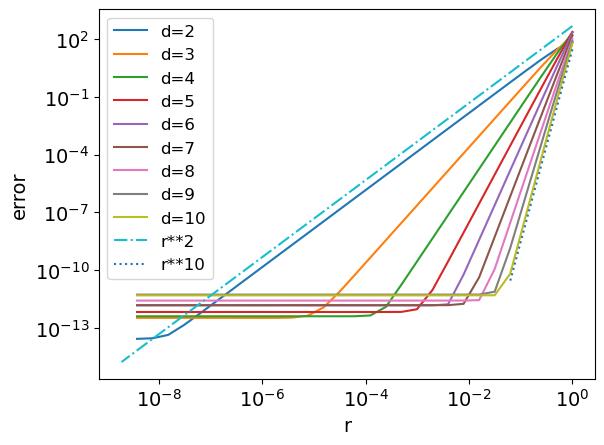}}
\caption{Maximum error over all elements on the pressure component $p$ between the exact solution $p_{ex}$ and its quasi-Trefftz approximation, evaluated for decreasing distances $r$ from the element center, Test Case~1. Computed for polynomial degree $\degp=10$ with the four quasi-Trefftz constructions (EXPL2, ALGE2, EXPL1, ALGE1).}\label{fig:approxv}
\end{figure}

Following the proof of Proposition \ref{prop:BAP},
to approximate the function $(p_{ex}, \bv_{ex})$ by a linear combination of quasi-Trefftz basis functions one can simply match the Taylor expansion of the exact solution by the Taylor expansion of a linear combination of quasi-Trefftz basis functions. In practice, given a quasi-Trefftz basis, one can build a matrix $\mathbf \phi_\degp$ of size 
$$
    \frac{(\degp+1)(\degp+2)+ 2\degp(\degp+1)}2\times (2\degp +1)
    =
    \frac{(\degp+1)(3\degp+2)}2\times (2\degp +1)
$$ 
whose columns each correspond to a basis function and whose rows each correspond to one Taylor expansion coefficients of the basis functions, understood as follows:
\begin{itemize}
    \item the Taylor expansion is taken component-wise,
    \item the expansion of the pressure component is taken up to degree $\degp$,
    \item the expansion of the velocity components is taken up to degree $\degp-1$.
\end{itemize}
The corresponding rectangular system can then be solved directly using a QR-decomposition.
Once the coefficients of the linear combination are determined, the approximation error is evaluated as follows. For the pressure component $p$, we compute the difference between the linear combination and the exact solution $p_{ex}$, at points located at decreasing distances $r$ from the element center. For each element, we record the maximum error over these points and over all basis functions, and then retain the largest value across all elements. The same procedure is applied to the velocity components $v_x$ and $v_y$, for which we take the worst of the two components. The resulting errors are reported in Figures~\ref{fig:approxp} and~\ref{fig:approxv}. 

We observe the convergence behavior expected from the theoretical result in Proposition \ref{prop:BAP}: for the pressure component $p$, the maximum error decreases with a slope of $\degp+1$ until reaching a plateau near machine precision; for the velocity components $(v_x,v_y)$ the slope corresponds to order $\degp$. 
It is worth noting that for the explicit methods, the plateau remains roughly at the same level across all degrees, whereas for the SVD-based methods the plateau tends to increase slightly with higher degrees. This likely reflects a worse conditioning of the matrix containing the Taylor coefficients of the basis functions in the algebraic construction.

\section{Beyond Helmholtz and polynomial bases}
\label{sec:GPWnHelmcomment}

The present work starts from a first-order formulation of the Helmholtz equation, as a step in the development of quasi-Trefftz methods for first-order differential systems. 
Next, the authors are particularly interested in tackling the convected Helmholtz equation for aeroacoustic applications. 
In this case, one preliminary challenge consists in choosing a first-order formulation of the PDE.
It is expected that both first-order equations will have variable coefficients, unlike for the Helmholtz case.

In \cite{LMIGpol}, a general framework developed to study polynomial quasi-Trefftz spaces for scalar PDEs was presented.
In \cite{GPWgen}, this general framework was extended to GPW quasi-Trefftz spaces for scalar PDEs.
Similarly, it would be interesting to extend the present work beyond the polynomial ansatz to a GPW ansatz. In this case, the quasi-Trefftz space cannot be characterized as the kernel of a linear operator. However, it is expected that, as for the scalar case, it would be possible to split the non-linear system defining a quasi-Trefftz property for the system case into a hierarchy of linear subsystems according to the graded structure of polynomial spaces.

\bibliographystyle{alpha}
\bibliography{biblio}

\end{document}